\documentclass[10pt]{article}

\usepackage{amsthm}
\usepackage{amsmath}
\usepackage{amsfonts}
\usepackage{amssymb}
\usepackage{mathrsfs}
\usepackage{empheq}
\usepackage{stmaryrd}
\usepackage{dsfont}
\usepackage{enumerate}
\usepackage{cancel}
\usepackage{comment}
\usepackage{pgfplots,soul}
\usepackage[applemac]{inputenc}
\usepackage{hyperref}

\usepackage{xcolor}

\newcommand{\N}{\mathbb{N}}

\newcommand{\R}{\mathbb{R}}
\newcommand{\C}{\mathbb{C}}

\newcommand{\vH}{{\cal H}}

\newcommand{\Dr}{\mathscr{D}}

\newcommand{\Lr}{\mathscr{L}}
\newcommand{\vPr}{\mathscr{Pr}}
\newcommand{\vphi}{\varphi}
\newcommand{\eps}{\varepsilon}

\newcommand{\dsp}{\displaystyle}
\newcommand{\ovl}{\overline}
\newcommand{\udl}{\underline}
\newcommand{\vlim}{\lim\limits}

\newcommand{\vsup}{\sup\limits}

\newcommand{\vint}{\int\limits}

\newcommand{\inj}{\hookrightarrow}
\newcommand{\tends}{\longrightarrow}
\newcommand{\weak}{\rightharpoonup}
\newcommand{\wt}{\widetilde}

\newcommand{\loc}{\mathrm{loc}}

\renewcommand{\b}{\mathrm{b}}
\newcommand{\co}{\mathrm{c}}
\renewcommand{\d}{\mathrm{d}}

\newcommand{\GN}{\mathrm{GN}}
\newcommand{\vi}{\mathrm{i}}
\newcommand{\w}{{\textsl w}}
\renewcommand{\le}{\leqslant}
\renewcommand{\ge}{\geqslant}
\renewcommand{\Re}{\mathrm{Re}}
\renewcommand{\Im}{\mathrm{Im}}
\newcommand{\bs}{\boldsymbol}
\newcommand{\p}{\prime}

\newcommand{\eqdef}{\stackrel{\mathrm{def}}{=}}
\DeclareMathOperator{\supp}{supp}

\DeclareMathOperator{\Arg}{Arg}

\numberwithin{equation}{section}

\newtheorem{thm}{Theorem}[section]
\newtheorem{prop}[thm]{Proposition}
\newtheorem{cor}[thm]{Corollary}
\newtheorem{lem}[thm]{Lemma}

\theoremstyle{definition}
\newtheorem{rmk}[thm]{Remark}
\newtheorem{defi}[thm]{Definition}

\newtheorem{assum}[thm]{Assumption}

\newenvironment{proof*}{\noindent{\bf Proof.}}{\qed}
\newenvironment{vproof}[1]{\noindent{\bf Proof #1}}{\qed}

\title{\huge \sc Finite time extinction for a class of damped Schr\"odinger equations with a singular saturated nonlinearity}
\author{\sc Pascal Bégout$^*$ and Jes\'us Ildefonso D{\'{\i}}az$^\dagger$}
\date{}

\begin{document}

\maketitle

\begin{gather*}
\begin{array}{cc}
                         ^*\mbox{Institut de Mathématiques de Toulouse}	&	\;^\dagger\mbox{Instituto de Matem\'atica Interdisciplinar}		\\
                                          \mbox{Université Toulouse I Capitole}	&	\mbox{Universidad Complutense de Madrid}				\\
                                             \mbox{1, Esplanade de l’Université}	&	\mbox{Plaza de las Ciencias, 3}						\\
                                  \mbox{31080 Toulouse Cedex 6, FRANCE}	&	\mbox{28040 Madrid, SPAIN}							\\
\bigskip \\
\mbox{
{\footnotesize E-mail\:: \href{mailto:Pascal.Begout@math.cnrs.fr}{\udl{\texttt{Pascal.Begout@math.cnrs.fr}}}}}
&
\mbox{
{\footnotesize E-mail\:: \href{mailto:jidiaz@ucm.es}{\udl{\texttt{jidiaz@ucm.es}}}}
}
\end{array}
\end{gather*}

\begin{abstract}
We present some sharper finite extinction time results for solutions of a class of damped nonlinear Schr\"{o}dinger equations when the nonlinear damping term corresponds to the limit cases of some ``saturating non-Kerr law'' $F(|u|^2)u=\frac{a}{\eps+(|u|^2)^\alpha}u,$ with $a\in\C,$ $\eps\ge0,$ $2\alpha=(1-m)$ and $m\in[0,1).$ To carry out the improvement of previous results in the literature we present in this paper a careful revision of the existence and regularity of weak solutions under very general assumptions on the data. We prove that the problem can be solved in the very general framework of the maximal monotone operators theory, even under a lack of regularity of the damping term. This allows us to consider, among other things, the singular case $m=0.$ We replace the above approximation of the damping term by a different one which keeps the monotonicity for any $\eps\ge0$. We prove that, when $m=0,$ the finite extinction time of the solution arises for merely bounded right hand side data $f(t,x).$ This is specially useful in the applications in which the Schr\"{o}dinger equation is coupled with some other functions satisfying some additional equations.
\end{abstract}

{\let\thefootnote\relax\footnotetext{2020 Mathematics Subject Classification: 35Q55 (35A01, 35A02, 35B40, 35D30, 35D35)}}
{\let\thefootnote\relax\footnotetext{Keywords: Damped Schr\"{o}dinger equation, Finite time extinction, Maximal monotone operators, Existence and regularity of weak solutions, Asymptotic behavior}}

\tableofcontents

\baselineskip .6cm

\section{Introduction}
\label{introduction}

This paper deals mainly with the asymptotic behavior, as $t\tends\infty,$ of solutions of the damped nonlinear Schr\"{o}dinger equation
\begin{empheq}[left=\empheqlbrace]{align}
	\label{nls}
	\vi\frac{\partial u}{\partial t}+\Delta u+V(x)u+a|u|^{-(1-m)}u=f(t,x),	&	\text{ in } (0,\infty)\times\Omega,				\\
	\label{nlsb}
	u_{|\partial\Omega}=0,									&	\text{ on } (0,\infty)\times\partial\Omega,	\dfrac{}{}	\\
	\label{u0}
	u(0)= u_0,												&	\text{ in } \Omega,
\end{empheq}
where $0\le m\le1,$ $a\in\C,$ $\Omega\subseteq\R^N$ (with $|\Omega|<\infty,$ if $m=0),$ $f\in L^1_\loc\big([0,\infty);L^2(\Omega)\big),$ $V\in L^1_\loc(\Omega;\R)$ and $u_0\in L^2(\Omega).$ More precisely, we will improve previous results in the literature (\cite{MR3306821}, \cite{MR4053613}, \cite{MR4098330}) showing that the main assumption 
\begin{gather}
\label{AssumpExtinction}
0\le m<1
\end{gather}
implies the finite extinction time phenomenon $(u(t)\equiv 0$ on $\Omega$ for any $t\ge T_\star,$ for some finite $T_\star>0)$ representing, clearly, the most opposite property to the famous Max Born result on the \textit{conservation of the mass}
\begin{gather*}
\|u(t)\|_{L^{2}(\Omega)}=\|u(s)\|_{L^{2}(\Omega)}, \text{ for any } t\ge s\ge 0,
\end{gather*}
which arises in the context of the applications of the linear Schr\"{o}dinger equation in Quantum Mechanics. Notice that this kind of non linear term can be understood as a special ``saturating non-Kerr law'' which arises in several applications (see, e.g. \cite{MR2000f:35139}, \cite{ak}, \cite{MR3308230}, \cite{MR1870804}, \cite{MR3470743} and their references) in which the following general nonlinear expression arises in the equation 
\begin{gather}
\label{approxF}
F(|u|^2)u=\frac{a}{\eps+(|u|^2)^\alpha}u,
\end{gather}
with $\eps\ge0$ and $2\alpha=(1-m).$ The assumption $m\in (0,1)$ corresponds to $\alpha\in(0,1/2)$ and the case $m=0$ corresponds to $\alpha=1/2.$ The consideration of the limit case, $\eps=0$ (assumed in this paper) allows to know the limit behavior of solutions for other weakly saturated cases in which $\eps>0.$ When $\eps=0$ the saturating term becomes singular at $u=0.$ We send the reader to the papers \cite{MR3306821}, \cite{MR4053613}, \cite{MR4098330} for many other information on the modeling and related results concerning problem \eqref{nls}--\eqref{u0}.
\medskip \\
It was already shown in the above mentioned works (\cite{MR3306821}, \cite{MR4053613}, \cite{MR4098330}) that the mere assumption \eqref{AssumpExtinction} is not enough to get to such a global conclusion and some other ``additional conditions'' are required. The improvements presented in this paper deal mainly with such type of ``additional conditions''. Some of them could be understood as ``technical conditions'' but they are of not of minor relevance since they require even an important revision already of the notion of solution of the problem. Thus, curiously enough, in some cases a proof about the asymptotic behavior requires to improve the basic framework of the existence and uniqueness of solutions. To carry out the improvement of previous results in the literature we present in this paper a careful revision of the existence and regularity of weak solutions under very general assumptions on the data. We prove that problem \eqref{nls}--\eqref{u0} can be solved in the very general framework of the maximal monotone operators theory on $L^2(\Omega),$ even under a lack of regularity of the damping term. This allows us to consider, among other things, the singular case $m=0.$ We replace the above approximation of the damping term~\eqref{approxF} by a different one:
\begin{gather*}
g_\eps^m(u)=(|u|^2+\eps)^{-\frac{1-m}2}u,
\end{gather*}
which keeps the monotonicity for any $\eps\ge0.$
\medskip \\
The motivation to include in the equation a given data $f(t,x)$ in the right hand side of the equation comes from the fact that very often the solution $u(t,x)$ of the Schr\"{o}dinger equation is coupled with other unknown term $v(t,x)$ satisfying, perhaps, a different PDE (of the type of the Maxwell equation, Poisson equation, conservation laws equation, etc.). Under suitable conditions (see, e.g., the energy methods applied to some coupled systems in \cite{MR2002i:35001}), it is possible that the coupled unknown $v(t,x)$ also presents a finite time extinction (see, e.g., \cite{MR2914580} for the case of a nonlinear Maxwell system) and this is the reason why we will assume in some of our results that the given data $f(t,x)$ satisfy such a property. 
\medskip \\
In this paper we will extend the formulation used in the previously mentioned papers to the case in which there is a linear potential term $Vu$ (in the philosophy of the Gross-Pitaevski models) in the equation and, which is perhaps less considered in the former literature, the limit case $m=0.$ We will understand the associated nonlinear operator as multivalued (see Definition~\ref{defsol}, Part~\ref{defsol3} below) and we will prove a curious result which was not noticed in (\cite{MR3306821}) where the case $m=0$ was also considered for a special formulation of problem \eqref{nls}--\eqref{u0} and for dimensions $N\le 2:$ the extinction time phenomenon holds in the larger class of data $f(t,x)$ for which we replace the condition $f(t)=0,$ a.e. $t\ge T_0,$ by the conditions
\begin{gather*}
f\in L^\infty\big((T_0,\infty)\times\Omega\big)
\end{gather*}
and 
\begin{gather}
\label{Hypofsmall}
\|f\|_{L^\infty((T_0,\infty)\times\Omega)}<\Im(a).
\end{gather}
In particular, when $f(t,x)$ represents a function of other possible coupled unknown $v(t,x),$ as mentioned before, this new condition is much more general than the assumption that $v(t,x)$ also presents a finite extinction time. We mention that although some related abstract results are available in the literature (see \cite{MR0420345} and \cite{DiazAcademia}) they can not be applied to the framework of problem \eqref{nls}--\eqref{u0}: see also this kind of property in the context of multivalued quasilinear parabolic equation (\cite{MR2466410}). Concerning multivalued hyperbolic wave equations, the phenomenon is associated to the presence of a Coulomb friction term in the equation (see, e.g., \cite{MR657297}, \cite{DiazMillot}, \cite{MR2185211} and \cite{MR2371117}) but usually $f(t,x)\equiv 0$ in this type of problems. See also the control point of view for some Maxwell class of scattering passive systems in \cite{SinghTucsnak}. The proof of our result (see Theorems~\ref{thmextN1} and \ref{thmextH2} below) is quite simple and avoids the application of any abstract result.
\medskip \\
Although very precise statements will be presented later, we point out that our Theorem~\ref{thmstrongaH1} below allows the consideration of data satisfying merely $f\in L^1(H^1_0)$ (and not necessarily $f\in W^{1,1}(H^1_0)$ as assumed in \cite{MR4053613} and \cite{MR4098330}). Another improvement of a ``technical nature'' is that in Theorems \ref{thmextN1e}, \ref{thmextH2} and \ref{thmextH2e}, we do not need to assume that $\Omega$ is a bounded
regular set if $\Omega\neq\R^N.$ Moreover, in the special case of $\Omega$ a half-space we show that $u,\Delta u\in L^2(\Omega)$ implies that $u\in H^2(\Omega),$ which is used in \ref{rmkthmtdH22} of Remark~\ref{rmkthmtdH2}.
\medskip \\
A different additional contribution, with respect to the previous papers (\cite{MR4053613} and \cite{MR4098330}) is that when we are not able to prove the finite extinction time at least we obtain some decay estimates as $t\tends\infty.$ For instance, we prove some cases in which $\vlim_{t\to\infty}\|\nabla u(t)\|_{L^2(\Omega)}=0$ (see Theorem~\ref{thmtdH2}).
\medskip \\
This paper is organized as follows. In Section~\ref{exiuni}, we state the main results about existence, uniqueness and boundness of solutions of~\eqref{nls}--\eqref{u0} (Theorem~\ref{thmweak}, \ref{thmstrongaH1}, \ref{thmstrongH1} and \ref{thmstrongH2}). In Section~\ref{finite}, we present the statements of our results about the finite time extinction property (Theorems~\ref{thmextN1}, \ref{thmextN1e}, \ref{thmextH2} and \ref{thmextH2e}) and on the asymptotic behavior (Theorems~\ref{thm0w}, \ref{thmrtdH1}, \ref{thmtdH1}, \ref{thmrtdH2} and \ref{thmtdH2}). Their respective proofs are also structured in different sections. In Section~\ref{functionals}, we give some a priori estimates about the terms $Vu$ and $|u|^{-(1-m)}u$ arising in the equation~\eqref{nls}, but, the more important part of the proofs is based on Section~\ref{maxmon} in which we will prove that it is possible to apply the theory of nonlinear maximal monotone operators on $L^2(\Omega)$ to the case of equation~\eqref{nls}. To this end we give some monotonicity results which are slight generalizations of the previous ones due to Liskevich and Perel$^\p$muter~\cite{MR1224619} and Hayashi~\cite{MR3802567}. Some additional properties and the proofs of the existence, uniqueness and boundedness of solutions are collected in Section~\ref{proofexi}. The paper ends with Section~\ref{proofext} with the proofs of
the statements on the finite extinction time and on the asymptotic behavior presented in Section~\ref{finite}.
\medskip \\
To end this Introduction, we collect here some notations which will be used along with this paper. For $t\in\R,$ $t_+=\max\{t,0\}$ is the positive part of $t.$ For $z\in\C,$ $\ovl z$ is the conjugate of $z,$ $\Re(z)$ is its real part and $\Im(z)$ is its imaginary part. Unless if specified, all functions are complex-valued $(H^1(\Omega)\eqdef H^1(\Omega;\C),$ etc) and all the vector spaces are considered over the field $\R.$ For $1\le p\le\infty,$ $p^\prime$ is the conjugate of $p$ defined by $\frac{1}{p}+\frac{1}{p^\prime}=1.$ For a (real)  Banach space $X,$ we denote by $X^\star\eqdef\Lr(X;\R)$ its topological dual and by $\langle\: . \; , \: . \:\rangle_{X^\star,X}\in\R$ the $X^\star-X$ duality product. In particular, for any $T\in L^{p^\prime}(\Omega)$ and $u\in L^p(\Omega)$ with $1\le p<\infty,$ $\langle T,u\rangle_{L^{p^\prime}(\Omega),L^p(\Omega)}=\Re\int_{\Omega}T(x)\ovl{u(x)}\d x.$ The scalar product in $L^2(\Omega)$ between two functions $u,v$ is, $(u,v)_{L^2(\Omega)}=\Re\int_{\Omega}u(x)\ovl{v(x)}\d x.$ For a Banach space $X$ and $p\in(0,\infty],$ $u\in L^p_\loc\big([0,\infty);X\big)$ means that $u\in L^p_\loc\big((0,\infty);X\big)$ and for any $T>0,$ $u_{|(0,T)}\in L^p\big((0,T);X\big).$ In the same way, we will use the notation $u\in W^{1,p}_\loc\big([0,\infty);X\big).$ If $p\in(0,\infty]$ then $L^\frac{p}r(\Omega)=L^\infty(\Omega)$ and $W^{1,\frac{p}r}(\Omega)=W^{1,\infty}(\Omega)$ if $r=0,$ and $L^0(\Omega)$ is the space of measurable functions $u:\Omega\tends\C$ such that $|u|<\infty,$ almost eveywhere in $\Omega.$ As usual, we denote by $C$ auxiliary positive constants, and sometimes, for positive parameters $a_1,\ldots,a_n,$ write as $C(a_1,\ldots,a_n)$ to indicate that the constant $C$ depends only on $a_1,\ldots,a_n$ and that dependence is continuous (we will use this convention for constants which are not denoted merely by ``$C$'').

\section{Main results}
\label{exiuni}

For $m\in[0,1],$ let us introduce the following sets of complex numbers:
\begin{align}
\label{Cm}
& C(m)=\Big\{z\in\C; \; \Im(z)>0 \text{ and } 2\sqrt m\Im(z)\ge(1-m)|\Re(z)|\Big\},			\\
\label{Dm}
& D(m)=\Big\{z\in\C; \; \Im(z)>0 \text{ and } 2\sqrt m\Im(z)=(1-m)\Re(z)\Big\}, \; 0<m<1,	\\
\label{Cmi}
& C_{\mathrm{int}}(m)=C(m)\setminus D(m), \; 0<m<1.
\end{align}
In the particular cases $m=0$ and $m=1,$ the set $C(m)$ becomes,
\begin{align*}
& C(0)=\Big\{z\in\C; \; \Re(z)=0 \text{ and } \Im(z)>0\Big\},	\\
& C(1)=\Big\{z\in\C; \; \Im(z)>0\Big\}.
\end{align*}

\noindent
Our main assumptions concerning the existence of the solutions are the following:
\begin{assum}
\label{ass1}
We assume the following.
\begin{gather}
\label{m}
0\le m\le1,										\\
\label{O}
\Omega \text{ is any nonempty open subset of } \R^N,	\\
\label{m0}
|\Omega|<\infty, \text{ if } m=0,	\\
\begin{cases}
\label{a}
a\in C(m),					&	\text{if } m\in\{0,1\},		\medskip \\
a\in C_{\mathrm{int}}(m),		&	\text{if } 0<m<1.			\medskip
\end{cases}
\\
\label{V}
V\in L^\infty(\Omega;\R)+L^{p_V}(\Omega;\R),
\end{gather}
where,
\begin{gather}
\label{pV}
p_V=
\begin{cases}
2,							&	\text{if } N=1,	\\
2+\beta, \text{ for some } \beta>0,	&	\text{if } N=2,	\\
N,							&	\text{if } N\ge3.
\end{cases}
\end{gather}
\medskip
\end{assum}

\noindent
Here and after, we shall always identify $L^2(\Omega)$ with its topological dual. Let us recall some important results of Functional Analysis. Let $E$ and $F$ be locally convex Hausdorff topological vector spaces. If $E\overset{e}{\inj}F$ with dense embedding then $F^\star\overset{e^\star}{\inj}E^\star,$ where $e^\star$ is the transpose of $e:$
\begin{gather*}
\forall L\in F^\star, \; \forall x\in E, \; \langle e^\star(L),x\rangle_{E^\star,E}=\langle L,e(x)\rangle_{F^\star,F}.
\end{gather*}
If, furthermore, $E$ is reflexive then the embedding $F^\star\overset{e^\star}{\inj}E^\star$ is dense. Often, $e$ is the identity function, so that $e^\star$ is nothing else but the restriction to $E$ of continuous linear forms on $F.$ For more details, see Trèves~\cite[Corollary~5, p.188; Corollary, p.199; Theorem~18.1, p.184]{MR2296978}. Let $A_1$ and $A_2$ two Banach spaces be such that $A_1,A_2\subset\vH$ for some Hausdorff topological vector space $\vH.$ Then $A_1\cap A_2$ and $A_1+A_2$ are Banach spaces where,
\begin{gather*}
\|a\|_{A_1\cap A_2}=\max\big\{\|a\|_{A_1},\|a\|_{A_2}\big\}	\; \text{ and } \;
\|a\|_{A_1+A_2}=\inf_{\left\{\substack{a=a_1+a_2 \hfill \\ (a_1,a_2)\in A_1\times A_2}\right.}\Big(\|a_1\|_{A_1}+\|a_2\|_{A_2}\Big).
\end{gather*}
If, in addition, $A_1\cap A_2$ is dense in both $A_1$ and $A_2$ then,
\begin{gather*}
\big(A_1\cap A_2\big)^\star=A_1^\star+A_2^\star \; \text{ and } \; \big(A_1+A_2\big)^\star=A_1^\star\cap A_2^\star.
\end{gather*}
See, for instance, Bergh and Löfström~\cite{MR0482275} (Lemma~2.3.1 and Theorem~2.7.1). We will often apply these results in the following cases. Let $0\le m\le1,$ let $X=H\cap L^{m+1}(\Omega),$ where $H=L^2(\Omega)$ or $H=H^1_0(\Omega),$ and let $Y$ be a Banach space such that $Y\inj L^p(\Omega)$ with dense embedding, for some $p\in[1,\infty).$ We then have,
\begin{gather}
\label{defX*}
X^\star=H^\star+L^\frac{m+1}m(\Omega),								\\
\label{X}
\Dr(\Omega)\inj X\inj L^{m+1}(\Omega) \text{ with both dense embeddings,}	\\
\label{X*}
L^\frac{m+1}m(\Omega)\inj X^\star\inj\Dr^\p(\Omega),					\\
\label{dualg}
\langle u,v\rangle_{Y^\star,Y}=\langle u,v\rangle_{L^{p^\p}(\Omega),L^p(\Omega)}
=\Re\vint_{\Omega} u(x)\ovl{v(x)}\d x,
\end{gather}
for any $u\in L^{p^\p}(\Omega)$ and $v\in Y.$ If $1<q<\infty$ and $p=2$ then by~Bégout and D\'{\i}az~\cite[Lemma~A.4]{MR4053613},
\begin{gather}
\label{Xcon}
L^q_\loc\big([0,\infty);Y\big)\cap W^{1,q^\p}_\loc\big([0,\infty);Y^\star\big)\inj C\big([0,\infty);L^2(\Omega)\big).
\end{gather}
By reflexivity of $\Dr(\Omega),$ the emdeddings $X^\star\inj\Dr^\p(\Omega)$ and $L^\frac{m+1}m(\Omega)\inj\Dr^\p(\Omega)$ are always dense. If $0<m\le1$ or if $|\Omega|<\infty$ then $X$ is reflexive and the embedding $L^\frac{m+1}m(\Omega)\inj X^\star$ is dense.
\medskip \\
We recall the definition of solution (\cite{MR4098330,MR4053613}), with a slight modification for $m=0,$ since it is not treated in \cite{MR4098330,MR4053613}.

\begin{defi}
\label{defsol}
Assume~\eqref{m}, \eqref{O}, \eqref{V} and \eqref{pV}. Let $a\in\C,$ $f\in L^1_\loc\big([0,\infty);L^2(\Omega)\big)$ and $u_0\in L^2(\Omega).$ Let us consider the following assertions.
\begin{enumerate}
\item
\label{defsol1}
$u\in L^{m+1}_\loc\big([0,\infty);H^1_0(\Omega)\cap L^{m+1}(\Omega)\big)\cap W^{1,\frac{m+1}m}_\loc\big([0,\infty);H^\star+L^\frac{m+1}m(\Omega)\big).$
\item
\label{defsol2}
For almost every $t>0,$ $\Delta u(t)\in H^\star.$
\item
\label{defsol3}
\begin{enumerate}
\item
\label{defsol3a}
If $m>0$ then $u$ satisfies~\eqref{nls} in $\Dr^\p\big((0,\infty)\times\Omega\big).$
\item
\label{defsol3b}
If $m=0$ then there exists $U\in L^\infty\big((0,\infty)\times\Omega\big)$ such that $\|U\|_{L^\infty((0,\infty)\times\Omega)}\le1,$ $U(t,x)=\dfrac{u(t,x)}{|u(t,x)|},$ if $u(t,x)\neq0,$ and $u$ satisfies~\eqref{nls} in $\Dr^\p\big((0,\infty)\times\Omega\big),$ where the term $|u|^{-(1-m)}u$ is replaced with $U.$
\end{enumerate}
\item
\label{defsol4}
$u(0)=u_0.$
\end{enumerate}
We shall say that $u$ is a \textit{strong solution} if $u$ is an $H^2$-solution or an $H^1_0$-solution. We shall say that $u$ is an $H^2$-\textit{solution of} \eqref{nls}--\eqref{u0} \big(respectively, an $H^1_0$-\textit{solution of} \eqref{nls}--\eqref{u0}\big), if $u$ satisfies the Assertions~\ref{defsol1}--\ref{defsol4} with $H=L^2(\Omega)$ \big(respectively, with $H=H^1_0(\Omega)\big).$
\\
We shall say that $u$ is an $L^2$-\textit{solution} or a \textit{weak solution} of \eqref{nls}--\eqref{u0} if there exists a pair,
\begin{gather}
\label{fn}
(f_n,u_n)_{n\in\N}\subset L^1_\loc\big([0,\infty);L^2(\Omega)\big)\times C\big([0,\infty);L^2(\Omega)\big),
\end{gather}
such that for any $n\in\N,$ $u_n$ is an $H^2$-solution of \eqref{nls}--\eqref{u0} where the right hand side of \eqref{nls} is $f_n,$ with
\begin{gather}
\label{cv}
f_n\xrightarrow[n\to\infty]{L^1((0,T);L^2(\Omega))}f \; \text{ and } \; u_n\xrightarrow[n\to\infty]{C([0,T];L^2(\Omega))}u,
\end{gather}
for any $T>0.$
\end{defi}

\noindent
Throughout this paper, we shall use the following notations and conventions. Let $m\in[0,1].$ Since $\big||z|^{-(1-m)}z\big|=|z|^m,$ we extend by continuity at $z=0$ the map $z\longmapsto|z|^{-(1-m)}z$ by setting,
\begin{gather*}
|z|^{-(1-m)}z=0, \text{ if } m>0 \text{ and } z=0.
\end{gather*}
Let $\eps\ge0.$ For any $u\in L^0(\Omega)$ and almost every $x\in\Omega,$ we define
\begin{align*}
& g_\eps^m(u)(x)=(|u(x)|^2+\eps)^{-\frac{1-m}2}u(x),	\; m+\eps>0,	\\
& g_0^0(u)(x)=\frac{u(x)}{|u(x)|}, \; u(x)\neq0,					\\
& g(u)(x)=g_0^m(u)(x).
\end{align*}

\begin{rmk}
\label{rmkdefsol}
Let us clarify the Definition~\ref{defsol}. See also Bégout~\cite{MR4098330}
for more details in the case $m>0.$
\begin{enumerate}
\item
\label{rmkdefsol1}
If $u$ is any strong or weak solution then by \ref{defsol1} of Definition~\ref{defsol}, \eqref{Xcon}, \eqref{cv} and the embedding
\begin{gather*}
 W^{1,\infty}_\loc\big([0,\infty);H^\star+L^\infty(\Omega)\big)\inj C\big([0,\infty);H^\star+L^\infty(\Omega)\big),
\end{gather*}
we have,
\begin{align*}
&	u\in C\big([0,\infty);L^2(\Omega)\big), \text{ if } m>0 \text{ or if } u \text{ is a weak solution,}							\\
&	u\in C\big([0,\infty);H^{-1}(\Omega)+L^\infty(\Omega)\big), \text{ if } m=0 \text{ and if } u \text{ is an } H^1_0\text{-solution,}		\\
&	u\in C\big([0,\infty);L^2(\Omega)+L^\infty(\Omega)\big), \text{ if } m=0 \text{ and if } u \text{ is an } H^2\text{-solution,}
\end{align*}
and thus the Cauchy condition $u(0)=u_0$ makes sense in some functional space according to the above different cases. Assume $m=0$ and $|\Omega|<\infty$ (such as indicated in Assumption~\ref{ass1}). Then, it is obvious that if $u$ is an $H^2$-solution then $u\in C\big([0,\infty);L^2(\Omega)\big).$ We claim that if $u$ is an $ H^1_0$-solution with $m=0,$ $|\Omega|<\infty$ and $f\in L^{1+\eps}_\loc\big([0,\infty);H^{-1}(\Omega)\big),$ for some $\eps\in(0,1),$ then for $r=2\frac{\eps+1}{\eps+2}\in(1,2),$
\begin{gather}
\label{rmkdefsol11}
u\in L^r_\loc\big([0,\infty);H^1_0(\Omega)\big)\cap W^{1,r^\p}_\loc\big([0,\infty);H^{-1}(\Omega)\big)\inj C\big([0,\infty);L^2(\Omega)\big).
\end{gather}
Indeed, by \ref{defsol1} of Definition~\ref{defsol} and the inequality $\|u(t)\|_{L^2(\Omega)}^2\le\|u(t)\|_{H^{-1}(\Omega)}\|u(t)\|_{H^1_0(\Omega)},$ we have $u\in L^2_\loc\big([0,\infty);L^2(\Omega)\big).$ With the help of \eqref{lemVH-1} below and \eqref{nls}, we get that $\Delta u\in L^{1+\eps}_\loc\big([0,\infty);H^{-1}(\Omega)\big).$ Finally, using the inequality $\|\nabla u(t)\|_{L^2(\Omega)}^r\le\|\Delta u(t)\|_{H^{-1}(\Omega)}^\frac{r}2\|u(t)\|_{H^1_0(\Omega)}^\frac{r}2,$ with $r=2\frac{\eps+1}{\eps+2},$ and integrating in time, we obtain, by the Hölder inequality,
\begin{gather*}
\|\nabla u\|_{L^r(0,T);L^2(\Omega))}^2\le\|\Delta u\|_{L^{1+\eps}((0,T);H^{-1}(\Omega))}\|u\|_{L^1((0,T);H^1_0(\Omega))},
\end{gather*}
for any $T>0.$ Hence \eqref{rmkdefsol11} holds.
\item
\label{rmkdefsol2}
Any $H^2$-solution satisfies~\eqref{nls} in $L^1_\loc\big([0,\infty);L^2(\Omega)+L^\frac{m+1}m(\Omega)\big),$ and any $H^1_0$-solution satisfies~\eqref{nls} in $L^1_\loc\big([0,\infty);H^{-1}(\Omega)+L^\frac{m+1}m(\Omega)\big).$ Indeed, this is a direct consequence of Definition~\ref{defsol} and Lemmas~\ref{lemVL} and \ref{lemcong} below.
\item
\label{rmkdefsol3}
Notice that the boundary condition $u(t)_{|\partial\Omega}=0$ is implicitely included in the assumption $u(t)\in H^1_0(\Omega),$ for the strong solutions. For the weak solutions, this has to be understood in a generalized sense by using the limit of strong solutions.
\end{enumerate}
\medskip
\end{rmk}

\noindent
The way in which the weak solutions satisfy the equation~\eqref{nls} is explained in the following result:

\begin{prop}
\label{propsolL2}
Let Assumption~$\ref{ass1}$ be fulfilled and let $f\in L^1_\loc\big([0,\infty);L^2(\Omega)\big).$ If $u$ is a weak solution to \eqref{nls} then
\begin{gather}
\label{propsolL21}
u\in W^{1,1}_\loc\big([0,\infty);H^{-2}(\Omega)+L^\frac2m(\Omega)\big).
\end{gather}
In addition, $u$ solves~\eqref{nls} in $L^1_\loc\big([0,\infty);H^{-2}(\Omega)+L^\frac2m(\Omega)\big)$ and so in $\Dr^\p\big((0,\infty)\times\Omega\big).$
\end{prop}

\noindent
Concerning the uniqueness and continuous dependance with respect to the initial data of solutions, we have:

\begin{prop}[\textbf{Uniqueness and continuous dependance}]
\label{propdep}
Assume~\eqref{m}--\eqref{m0} and \eqref{V}--\eqref{pV}. Let $a\in C(m),$ let $f,\wt f\in L^1_\loc\big([0,\infty);L^2(\Omega)\big)$ and $X=H^1_0(\Omega)\cap L^{m+1}(\Omega).$ Finally, let
\begin{gather}
\label{propdephyp}
u,\wt u\in L^p_\loc\big([0,\infty);X\big)\cap W^{1,p^\p}_\loc\big([0,\infty);X^\star\big)\inj C\big([0,\infty);L^2(\Omega)\big),
\end{gather}
for some $1<p<\infty,$ be solutions in $\Dr^\p\big((0,\infty)\times\Omega\big)$ to,
\begin{gather*}
\vi u_t+\Delta u+Vu+a|u|^{-(1-m)}u=f,		\\
\vi\wt {u_t}+\Delta\wt u+V\wt u+a|\wt u|^{-(1-m)}\wt u=\wt f ,
\end{gather*}
respectively $($with the obvious modification, as in Definition~$\ref{defsol},$ if $m=0).$ Then,
\begin{gather}
\label{estthmweak}
\|u(t)-\wt u(t)\|_{L^2(\Omega)}\le\|u(s)-\wt u(s)\|_{L^2(\Omega)}+\vint_s^t\|f(\sigma)-\wt f(\sigma)\|_{L^2(\Omega)}\d\sigma,
\end{gather}
for any $t\ge s\ge0.$
\end{prop}

\begin{thm}[\textbf{Existence and uniqueness of $\bs{L^2}$-solutions}]
\label{thmweak}
Let Assumption~$\ref{ass1}$ be fulfilled and let $f\in L^1_\loc\big([0,\infty);L^2(\Omega)\big).$ Then for any $u_0\in L^2(\Omega),$ there exists a unique weak solution $u$ to \eqref{nls}--\eqref{u0}. In addition,
\begin{gather}
\label{Lm}
u\in L^{m+1}_\loc\big([0,\infty);L^{m+1}(\Omega)\big),	\\
\label{L2+}
\dfrac12\|u(t)\|_{L^2(\Omega)}^2+\Im(a)\dsp\vint_s^t\|u(\sigma)\|_{L^{m+1}(\Omega)}^{m+1}\d\sigma
			\le\dfrac12\|u(s)\|_{L^2(\Omega)}^2+\Im\dsp\iint\limits_{s\;\Omega}^{\text{}\;\;t}f(\sigma,x)\,\ovl{u(\sigma,x)}\,\d x\,\d\sigma,
\end{gather}
for any $t\ge s\ge0.$ If $|\Omega|<\infty$ or if $m=1$ then the inequality in~\eqref{L2+} is an equality. Finally, if $\wt u$ is a weak solution to \eqref{nls} with $\wt u(0)=\wt{u_0}\in L^2(\Omega)$ and $\wt f\in L^1_\loc([0,\infty);L^2(\Omega))$ instead of $f$ in \eqref{nls} then~\eqref{estthmweak} holds for any $t\ge s\ge0.$
\end{thm}

\begin{thm}[\textbf{Additional regularity in $\bs{H^1_0}$}]
\label{thmstrongaH1}
Let Assumption~$\ref{ass1}$ be fulfilled. Assume that each component of the vector $\nabla V$ satisfies the regularity~\eqref{V} and let $f\in L^1_\loc\big([0,\infty);H^1_0(\Omega)\big).$ Then for any $u_0\in H^1_0(\Omega),$ the weak solution $u$ satisfies, additionally, that
\begin{gather}
\begin{cases}
\label{thmstrongaH11}
u\in C\big([0,\infty);L^2(\Omega)\big)\cap L^\infty_\loc\big([0,\infty);H^1_0(\Omega)\big)\cap L^{m+1}_\loc\big([0,\infty);L^{m+1}(\Omega)\big),	\medskip \\
u\in W^{1,1}_\loc\big([0,\infty);H^{-1}(\Omega)+L^\frac{m+1}m(\Omega)\big),
\end{cases}
\end{gather}
and $u$ satisfies \eqref{nls} in $ L^1_\loc\big([0,\infty);H^{-1}(\Omega)+L^\frac{m+1}m(\Omega)\big).$ In addition, $u$ verifies,
\begin{gather}
\label{thmstrongaH12}
\|u(t)\|_{H^1_0(\Omega)}\le\left(\|u(s)\|_{H^1_0(\Omega)}+\vint_s^t\|f(\sigma)\|_{H^1_0(\Omega)}\d\sigma\right)e^{C\|\nabla V\|_{L^\infty+L^{p_V}}(t-s)},
\end{gather}
for almost every $t>s>0,$ where $C=C(N)$ $(C=C(\beta),$ if $N=2).$
\end{thm}

\begin{rmk}
\label{rmkthmstrongaH1}
Below are some comments about Theorems~\ref{thmweak} and \ref{thmstrongaH1}.
\begin{enumerate}
\item
\label{rmkthmstrongaH11}
Let Assumption~$\ref{ass1}$ be fulfilled and let $u$ be a weak solution. If $f\in L^1\big((0,\infty);L^2(\Omega)\big)$ then,
\begin{gather}
\label{rmkthmweak1}
u\in C_\b\big([0,\infty);L^2(\Omega)\big)\cap L^\frac{p(1-m)}{2-p}\big((0,\infty);L^p(\Omega)\big),
\end{gather}
for any $p\in[m+1,2].$ Here, by $C_\b$ we mean $C\cap L^\infty.$ If, in addition, $(\vphi_n)_{n\in\N}\subset L^2(\Omega),$ $(f_n)_{n\in\N}\subset L^1\big((0,\infty);L^2(\Omega)\big)$ and,
\begin{gather*}
\vphi_n\xrightarrow[n\to\infty]{L^2(\Omega)}u_0 \; \text{ and } \; f_n\xrightarrow[n\to\infty]{L^1((0,\infty);L^2(\Omega))}f,
\end{gather*}
then for any $p\in(m+1,2),$
\begin{gather*}
u_n\xrightarrow[n\to\infty]{C_\b([0,\infty);L^2(\Omega))\cap L^\frac{p(1-m)}{2-p}((0,\infty);L^p(\Omega))}u,
\end{gather*}
where for each $n\in\N,$ $u_n$ is the weak solution to \eqref{nls} with $u_n(0)=\vphi_n$ and $f_n$ instead of $f.$ See Bégout~\cite[Remark~2.5]{MR4098330}
for more details.
\item
\label{rmkthmstrongaH12}
The solution obtained in Theorem~\ref{thmstrongaH1} could be called an \textit{almost} $H^1_0$\textit{-solution} since it verifies all the conditions of Definition~\ref{defsol}, except the property
\begin{gather}
\label{rmkthmstrongaH121}
u\in W^{1,\frac{m+1}m}_\loc\big([0,\infty);X^\star\big),
\end{gather}
which need not be satisfied, where $X^\star=H^{-1}(\Omega)+L^\frac{m+1}m(\Omega).$ In particular, we cannot apply Proposition~\ref{propdep} and, as a consequence, we do not know if the solution is unique in the class of functions satisfying~\eqref{thmstrongaH11}. Of course, it is unique in the class of weak solutions (Theorem~\ref{thmweak}). Finally, \eqref{rmkthmstrongaH121} may be obtained if we assume additionally $f\in L^\frac{m+1}m_\loc\big([0,\infty);X^\star\big),$ (see Theorem~\ref{thmstrongH1} below).
\item
\label{rmkthmstrongaH13}
The assumption on $\nabla V$ in Theorem~\ref{thmstrongaH1} (and Theorem~\ref{thmstrongH1} below) is needed to obtain \eqref{thmstrongaH12} and, thereby, the approximating sequence of the $H^2$-solutions bounded in $L^\infty_\loc\big([0,\infty);H^1_0(\Omega)\big).$ If $V$ is a constant function then we may obtain a better estimate as follows. We claim that,
\begin{gather}
\label{thmstrongaH12V0}
\|\nabla u(t)\|_{L^2(\Omega)}\le\|\nabla u(s)\|_{L^2(\Omega)}+\vint_s^t\|\nabla f(\sigma)\|_{L^2(\Omega)}\d\sigma,
\end{gather}
for almost every $t>s>0.$ Indeed, since the solution obtained in Theorem~\ref{thmstrongaH1} is a weak solution, by uniqueness of the weak solutions and by a time translation argument, it is sufficient to establish~\eqref{thmstrongaH12V0} for $s=0$ and the $H^2$-solutions. Taking the $L^2$-scalar product of \eqref{nls} with $-\vi\Delta u,$ it follows from~Bégout and D\'{\i}az~\cite[Lemma~A.5]{MR4053613}
and Lemma~\ref{lemenereg} below that for almost every $\sigma>0,$
\begin{gather*}
\frac12\frac{\d}{\d t}\|\nabla u(\sigma)\|_{L^2(\Omega)}^2\le\big(\nabla f(\sigma),\vi\nabla u(\sigma)\big)_{L^2(\Omega)}
\le\|\nabla f(\sigma)\|_{L^2(\Omega)}\|\nabla u(\sigma)\|_{L^2(\Omega)}.
\end{gather*}
The result then follows by integration. See the proof of Theorem~\ref{thmstrongaH1} for more details.
\end{enumerate}
\end{rmk}

\noindent
Below and after, we denote by $C_\w\big([0,\infty);H^1_0(\Omega)\big)$ the space of continuous functions from $[0,\infty)$ to $H^1_0(\Omega),$ where $H^1_0(\Omega)$ is endowed of the weak topology $\sigma\left(H^1_0(\Omega),H^{-1}(\Omega)\right).$

\begin{thm}[\textbf{Existence and uniqueness of $\bs{H^1_0}$-solutions}]
\label{thmstrongH1}
Let Assumption~$\ref{ass1}$ be fulfilled. Assume that each component of the vector $\nabla V$ satisfies~\eqref{V} and let
\begin{gather*}
f\in L^1_\loc\big([0,\infty);H^1_0(\Omega)\big)\cap L^\frac{m+1}m_\loc\big([0,\infty);H^{-1}(\Omega)+L^\frac{m+1}m(\Omega)\big).
\end{gather*}
Then for any $u_0\in H^1_0(\Omega),$ there exists a unique $H^1_0$-solution $u$ to \eqref{nls}--\eqref{u0}. Furthermore, $u$ is also a weak solution and satisfies the following properties.
\begin{enumerate}
\item
\label{thmstrongH11}
$u\in C_\w\big([0,\infty);H^1_0(\Omega)\big)$ and \eqref{thmstrongaH12} holds for any $t\ge s\ge0.$
\item
\label{thmstrongH12}
The map $t\longmapsto\|u(t)\|_{L^2(\Omega)}^2$ belongs to $W^{1,1}_\loc\big([0,\infty);\R\big)$ and we have,
\begin{gather}
\label{L2}
\frac12\frac{\d}{\d t}\|u(t)\|_{L^2(\Omega)}^2+\Im(a)\|u(t)\|_{L^{m+1}(\Omega)}^{m+1}=\Im\vint_{\Omega}f(t,x)\,\ovl{u(t,x)}\,\d x,
\end{gather}
for almost every $t>0.$
\end{enumerate}
\end{thm}

\begin{thm}[\textbf{Existence and uniqueness of $\bs{H^2}$-solutions}]
\label{thmstrongH2}
Let Assumption~$\ref{ass1}$ be fulfilled and $f\in W^{1,1}_\loc\big([0,\infty);L^2(\Omega)\big).$ Then for any $u_0\in H^1_0(\Omega)\cap L^{2m}(\Omega),$ with $\Delta u_0\in L^2(\Omega),$ there exists a unique $H^2$-solution $u$ to \eqref{nls}--\eqref{u0}. Furthermore, $u$ satisfies~\eqref{nls} in $L^\infty_\loc\big([0,\infty);L^2(\Omega)\big)$ as well as the following properties.
\begin{enumerate}
\item
\label{thmstrongH21}
$u\in C\big([0,\infty);H^1_0(\Omega)\big)\cap W^{1,\infty}_\loc\big([0,\infty);L^2(\Omega)\big)$ and, in addition, $u\in L^\infty_\loc\big([0,\infty);L^{2m}(\Omega)\big),$ if $m>0.$
\item
\label{thmstrongH22}
$\Delta u\in L^\infty_\loc\big([0,\infty);L^2(\Omega)\big)$ and,
\begin{empheq}[left=\empheqlbrace]{align}
\label{strongH21}
	&	\|u(t)-u(s)\|_{L^2(\Omega)}\le\|u_t\|_{L^\infty((s,t);L^2(\Omega))}|t-s|,			\dfrac{}{}			\\
\label{strongH22}
	&	\|\nabla u(t)-\nabla u(s)\|_{L^2(\Omega)}\le M|t-s|^\frac12,									\\
\label{strongH23}
	&	\left\|u_t\right\|_{L^\infty((0,t);L^2(\Omega))}
		\le\|\vi A_0^mu_0-f(0)\|_{L^2(\Omega)}+\int_0^t\|f^\p(\sigma)\|_{L^2(\Omega)}\d\sigma,
\end{empheq}
for any $t\ge s\ge0,$ where $M^2=2\|u_t\|_{L^\infty((s,t);L^2(\Omega))}\|\Delta u\|_{L^\infty((s,t);L^2(\Omega))}$ and $\vi A_0^mu_0=\Delta u_0+Vu_0+ag(u_0)$ $\big(\vi A_0^0u_0=\Delta u_0+Vu_0+aU_0,$ for some $U_0$ in the closed unit ball of $L^\infty(\Omega)$ with $U_0=\frac{u_0}{|u_0|},$ almost everywhere where $u_0\neq0,$ if $m=0\big).$
\item
\label{thmstrongH23}
The map $t\longmapsto\|u(t)\|_{L^2(\Omega)}^2$ belongs to $C^1\big([0,\infty);\R\big)$ and~\eqref{L2} holds for any $t\ge0.$
\item
\label{thmstrongH24}
If $f\in W^{1,1}\big((0,\infty);L^2(\Omega)\big)$ then we have,
\begin{align*}
	&	u\in C_\b\big([0,\infty);H^1_0(\Omega)\big)\cap W^{1,\infty}\big((0,\infty);L^2(\Omega)\big),	\\
	&	\Delta u\in L^\infty\big((0,\infty);L^2(\Omega)\big),									\\
	&	u\in L^\infty\big((0,\infty);L^{2m}(\Omega)\big), \text{ if } m>0.
\end{align*}
\end{enumerate}
\end{thm}

\begin{rmk}
\label{rmkf0}
Below are some comments about Theorem~\ref{thmstrongH2}.
\begin{enumerate}
\item
\label{rmkf01}
Since $f\in W^{1,1}_\loc\big([0,\infty);L^2(\Omega)\big)\inj C\big([0,\infty);L^2(\Omega)\big),$ estimate~\eqref{strongH23} with $f(0)$ makes sense.
\item
\label{rmkf02}
For any $p\in\left(2m,\frac{2N}{N-2}\right]$ $(p\in(2m,\infty)$ if $N=2,$ $p\in(2m,\infty]$ if $N=1),$
\begin{gather*}
u\in C^{0,\alpha}\big([0,\infty);L^p(\Omega)\big)	\;
\Big(u\in C^{0,\alpha}_\b\big([0,\infty);L^p(\Omega)\big), \text{ if } f\in W^{1,1}\big((0,\infty);L^2(\Omega)\big)\Big),
\end{gather*}
where $\alpha=\frac{2N-p(N-4)}{4p}$ if $p\ge2,$ $\alpha=\frac{p-2m}{p(1-m)}$ if $p\le2$ and $m>0,$ and $\alpha=1$ if $p\le2$ and $m=0.$ Indeed, if $p\ge2$ this comes from Properties~\ref{thmstrongH21} and \ref{thmstrongH22}, and Gagliardo-Nirenberg's inequality. If $m>0$ and $p\in(2m,2],$ this regularity comes from Hölder's inequality, Property~\ref{thmstrongH21} and \eqref{strongH21}. Finally, if $m=0$ and $p\in(0,2],$ this comes from \eqref{strongH21} and the embedding $L^2(\Omega)\inj L^p(\Omega),$ since $|\Omega|<\infty.$
\end{enumerate}
\end{rmk}

\begin{rmk}
\label{rmkDm}
The existence of the solutions of \eqref{nls}--\eqref{u0} for $a\in D(m)$ is not treated here and will be the subject of a future work. Note that if $|\Omega|<\infty$ and $V=0,$ this was done in~Bégout and D\'{\i}az~\cite{MR4053613}.
\end{rmk}

\section{Finite time extinction and asymptotic behavior}
\label{finite}

We will improve the result of~Bégout and D\'{\i}az~\cite{MR4053613}
by avoiding, among other things, some regularity and boundedness conditions on the spatial domain.
\medskip \\
For $N\in\N,$ let $\ell\in\{1,2\},$ $m\in[0,1)$ and
\begin{gather}
\label{delta}
\delta_\ell=\frac{(N+2\ell)-m(N-2\ell)}{4\ell}.
\end{gather}
Notice that if $N=\ell=1$ or if $N\le3$ then for any $m\in[0,1),$ $\delta_\ell\in\left(\frac12,1\right).$

\begin{assum}[\textbf{Case of the }$\bs{H^1_0}$\textbf{-solutions}]
\label{assN1}
Assumption~\ref{ass1} holds true with $0\le m<1$ and $V$ a constant function. Let $f\in L^1\big((0,\infty);H^1_0(\Omega)\big),$ let $u_0\in H^1_0(\Omega)$ and let $u$ be the unique $L^2$-solution to \eqref{nls}--\eqref{u0} given by Theorem~\ref{thmweak}. We assume that there exists a finite time $T_0\ge0$ such that
\begin{gather}
\label{f}
\begin{cases}
f\in L^\infty\big((T_0,\infty)\times\Omega\big) \text{ and } \|f\|_{L^\infty((T_0,\infty)\times\Omega)}<\Im(a),		&	\text{if } m=0,		\medskip \\
f(t)=0, \text{ for almost every } t>T_0,															&	\text{if } 0<m<1.
\end{cases}
\end{gather}
\end{assum}

\begin{assum}[\textbf{Case of the }$\bs{H^2}$\textbf{-solutions}]
\label{assH2}
Assumption~\ref{ass1} holds true with $0\le m<1.$ Let $f\in W^{1,1}\big((0,\infty);L^2(\Omega)\big),$ $u_0\in H^1_0(\Omega)\cap L^{2m}(\Omega)$ with $\Delta u_0\in L^2(\Omega)$ and let $u$ be the unique $H^2$-solution to \eqref{nls}--\eqref{u0} given by Theorem~\ref{thmstrongH2}. We assume that there exists a finite time $T_0\ge0$ such that $f$ satisfies~\eqref{f}.
\end{assum}

\section*{Asymptotic behavior of the $\bs{L^2}$-solutions}

\begin{thm}
\label{thm0w}
Let Assumption~$\ref{ass1}$ be fulfilled, let $f\in L^1\big((0,\infty);L^2(\Omega)\big),$ $u_0\in L^2(\Omega)$ and let $u$ be the unique weak solution to \eqref{nls}--\eqref{u0} given by Theorem~$\ref{thmweak}.$ Then,
\begin{gather*}
\vlim_{t\nearrow\infty}\|u(t)\|_{L^2(\Omega)}=0.
\end{gather*}
\end{thm}

\begin{rmk}
\label{rmkthm0w}
Let the hypotheses of Theorem~\ref{thm0w} be fulfilled with $m=0$ and $|\Omega|<\infty.$ By the embedding $L^2(\Omega)\inj L^p(\Omega),$ we have $\vlim_{t\nearrow\infty}\|u(t)\|_{L^p(\Omega)}=0,$ for any $p\in(0,2].$ Now, suppose $m=1$ and $f=0$ almost everywhere on $(T_0,\infty),$ for some $T_0\ge0.$ Then,
\begin{gather*}
\forall t\ge T_0, \; \|u(t)\|_{L^2(\Omega)}=\|u(T_0)\|_{L^2(\Omega)}e^{-\Im(a)(t-T_0)}.
\end{gather*}
Indeed, by~\eqref{estthmweak} and density, we may assume that $u$ is an $H^2$-solution. We then have by \eqref{L2},
\begin{gather*}
\forall t\ge T_0, \; \frac12\frac{\d}{\d t}\|u(t)\|_{L^2(\Omega)}^2+\Im(a)\|u(t)\|_{L^2(\Omega)}^2=0,
\end{gather*}
from which the result follows. We have a similar statement for the strong solutions when $m<1$ (Theorems~\ref{thmrtdH1} and \ref{thmrtdH2} below).
\end{rmk}

\section*{Finite time extinction and asymptotic behavior of the $\bs{H^1_0}$-solutions}

\begin{thm}[\textbf{Finite time extinction}]
\label{thmextN1}
Let Assumption~$\ref{assN1}$ be fulfilled with $N=1.$ Then,
\begin{gather}
\label{01}
\forall t\ge T_\star, \; \|u(t)\|_{L^2(\Omega)}=0,
\end{gather}
where,
\begin{gather}
\label{T*N1}
T_\star\le C\|u(T_0)\|_{L^2(\Omega)}^\frac{1-m}2\|\nabla u\|_{L^\infty((0,\infty);L^2(\Omega))}^\frac{1-m}2+T_0,
\end{gather}
for some $C=C(\Im(a),m)$ $(C=C(\Im(a)-\|f\|_{L^\infty((T_0,\infty)\times\Omega)}),$ if $m=0).$
\end{thm}

\begin{thm}[\textbf{Synchronized finite time extinction}]
\label{thmextN1e}
Let  Assumption~$\ref{assN1}$ be fulfilled with $N=1.$ Assume further that $f\in L^\frac{m+1}m\big((0,\infty);H^{-1}(\Omega)+L^\frac{m+1}m(\Omega)\big)$ so that $u$ is an $H^1_0$-solution. There exists $\eps_\star=\eps_\star(|a|,m)$ satisfying the following property. If
\begin{gather}
\label{thmextN1e1}
\begin{cases}
\|u_0\|_{L^2(\Omega)}^{2(1-\delta_1)}\le\eps_\star T_0,							\medskip \\
\|\nabla u_0\|_{L^2(\Omega)}+\|\nabla f\|_{L^1((0,\infty);L^2(\Omega))}\le\eps_\star,		\medskip \\
\|f(t)\|_{L^2(\Omega)}^2\le\eps_\star\big(T_0-t\big)_+^\frac{2\delta_1-1}{1-\delta_1},
\end{cases}
\end{gather}
for almost every $t>0,$ where $\delta_1$ is defined by~\eqref{delta}, then~\eqref{01} holds true with $T_\star=T_0.$
\end{thm}

\begin{thm}[\textbf{Time decay estimates}]
\label{thmrtdH1}
Let Assumption~$\ref{assN1}$ be fulfilled with $N\ge2.$ Then for any $t\ge T_0,$
\begin{gather}
\label{thmrtdH11}
\|u(t)\|_{L^2(\Omega)}\le\|u(T_0)\|_{L^2(\Omega)}e^{-C(t-T_0)},
\end{gather}
if $N=2,$ and
\begin{gather}
\label{thmrtdH12}
\|u(t)\|_{L^2(\Omega)}\le\dfrac{\|u(T_0)\|_{L^2(\Omega)}}
{\left(1+C\|u(T_0)\|_{L^2(\Omega)}^\frac{(1-m)(N-2)}2(t-T_0)\right)^\frac2{(1-m)(N-2)}},
\end{gather}
if $N\ge3,$ where $C=C(\|\nabla u\|_{L^\infty((0,\infty);L^2(\Omega))},\Im(a),N,m)$ $(C=C(\|\nabla u\|_{L^\infty((0,\infty);L^2(\Omega))},\Im(a)-\|f\|_{L^\infty((T_0,\infty)\times\Omega)},N),$ if $m=0).$
\end{thm}

\begin{thm}[\textbf{Time decay}]
\label{thmtdH1}
Let Assumption~$\ref{ass1}$ be fulfilled with $V$ a constant function. Let $f\in L^1\big((0,\infty);H^1_0(\Omega)\big),$ $u_0\in H^1_0(\Omega)$ and let $u$ be the unique weak solution given by Theorem~$\ref{thmweak}.$ Then,
\begin{gather}
\label{thmtdH11}
\vlim_{t\nearrow\infty}\|u(t)\|_{L^p(\Omega)}=0,
\end{gather}
for any $p\in\left[2,\frac{2N}{N-2}\right)$ $(p\in[2,\infty]$ if $N=1).$ If $m=0$ then \eqref{thmtdH11} is also true for any $p\in(0,2].$
\end{thm}

\section*{Finite time extinction and asymptotic behavior of the $\bs{H^2}$-solutions}

\begin{thm}[\textbf{Finite time extinction}]
\label{thmextH2}
Let Assumption~$\ref{assH2}$ be fulfilled with $N\le3.$ Then,
\begin{gather}
\label{02}
\forall t\ge T_\star, \; \|u(t)\|_{L^2(\Omega)}=0,
\end{gather}
where,
\begin{gather}
\label{T*H2}
T_\star\le C\|u(T_0)\|_{L^2(\Omega)}^\frac{(1-m)(4-N)}4\|\Delta u\|_{L^\infty((0,\infty);L^2(\Omega))}^\frac{N(1-m)}4+T_0,
\end{gather}
for some $C=C(\Im(a),N,m)$ $(C=C(\Im(a)-\|f\|_{L^\infty((T_0,\infty)\times\Omega)},N),$ if $m=0).$
\end{thm}

\begin{rmk}
\label{rmkUf0}
Assume $m=0.$ When $u(t,x)=0,$ we do not know exactly what is the term $U(t,x)$ in the equation~\eqref{nls} (remember Part~\ref{defsol3b} of Definition~\ref{defsol}), except in the following particular case. Assume that, for some $T_0\ge0,$ $f$ satisfies~\eqref{f}. Let $u$ be a solution as in Theorems~\ref{thmextN1} or \ref{thmextH2}. Then by~\eqref{01} or \eqref{02}, the equation~\eqref{nls} becomes,
\begin{gather*}
\vi\Im(a)\,U(t,x)=f(t,x),
\end{gather*}
for almost every $(t,x)\in(T_\star,\infty)\times\Omega.$
\end{rmk}

\noindent
For $0<m<1,$ let us define the quasi-norm $\|\:.\:\|_{m,\Omega}$ by,
\begin{gather}
\label{mn}
\|u\|_{m,\Omega}=\|u\|_{H^1_0(\Omega)}+\|u\|_{L^{2m}(\Omega)}+\|\Delta u\|_{L^2(\Omega)},
\end{gather}
for any $u\in H^1_0(\Omega)\cap L^{2m}(\Omega)$ with $\Delta u\in L^2(\Omega).$

\begin{thm}[\textbf{Synchronized finite time extinction}]
\label{thmextH2e}
Let  Assumption~$\ref{assH2}$ be fulfilled with $N\le3$ and $0<m<1.$ There exists $\eps_\star=\eps_\star(|a|,N,m)$ satisfying the following property. If
\begin{gather}
\label{thmextH2e1}
\begin{cases}
\|u_0\|_{L^2(\Omega)}^{2(1-\delta_2)}\le\eps_\star T_0,							\medskip \\
\|u_0\|_{m,\Omega}+\|f\|_{W^{1,1}((0,\infty);L^2(\Omega))}\le\eps_\star,				\medskip \\
\|f(t)\|_{L^2(\Omega)}^2\le\eps_\star\big(T_0-t\big)_+^\frac{2\delta_2-1}{1-\delta_2},
\end{cases}
\end{gather}
for almost every $t>0,$ where $\delta_2$ is defined by~\eqref{delta}, then~\eqref{02} holds true with $T_\star=T_0.$
\end{thm}

\begin{thm}[\textbf{Time decay estimates}]
\label{thmrtdH2}
Let Assumption~$\ref{assH2}$ be fulfilled with $N\ge4.$ Then for any $t\ge T_0,$
\begin{gather}
\label{thmrtdH21}
\|u(t)\|_{L^2(\Omega)}\le\|u(T_0)\|_{L^2(\Omega)}e^{-C(t-T_0)},
\end{gather}
if $N=4,$ and
\begin{gather}
\label{thmrtdH22}
\|u(t)\|_{L^2(\Omega)}\le\dfrac{\|u(T_0)\|_{L^2(\Omega)}}
{\left(1+C\|u(T_0)\|_{L^2(\Omega)}^\frac{(1-m)(N-4)}4(t-T_0)\right)^\frac4{(1-m)(N-4)}},
\end{gather}
if $N\ge5,$ where $C=C(\|\Delta u\|_{L^\infty((0,\infty);L^2(\Omega))},\Im(a),N,m)$ $(C=C(\|\Delta u\|_{L^\infty((0,\infty);L^2(\Omega))},\Im(a)-\|f\|_{L^\infty((T_0,\infty)\times\Omega)},N),$ if $m=0).$
\end{thm}

\begin{rmk}
\label{rmkfm0}
As mentioned at the introduction, the results of Theorems~\ref{thmextN1}, \ref{thmrtdH1}, \ref{thmextH2} and \ref{thmrtdH2} for $m=0$ can be applied to the case in which the Schrödinger equation is coupled with some other dynamic equation 
\begin{gather*}
\begin{cases}
\vi\dfrac{\partial u}{\partial t}+\Delta u+V(x)u+a\dfrac{u}{|u|}=g(v),	&	\medskip \\
\dfrac{\partial v}{\partial t}+B(v)=h(u,v).						&
\end{cases}
\end{gather*}
Then by taking $f(t,x)=g\big(v(t,x)\big),$ if we can prove, for instance, that $\|v(t)\|_{L^\infty(\Omega)}\xrightarrow{t\to\infty}0,$ and if $g$ is Lipschitz continuous with $g(0)=0,$ then the assumption $\|f\|_{L^\infty((T_0,\infty)\times\Omega)}<\Im(a)$ is satisfied, for $T_0>0$ large enough.
\end{rmk}

\begin{thm}[\textbf{Time decay}]
\label{thmtdH2}
Let Assumption~$\ref{ass1}$ be fulfilled. Let $f\in W^{1,1}\big((0,\infty);L^2(\Omega)\big),$ $u_0\in H^1_0(\Omega)\cap L^{2m}(\Omega)$ with $\Delta u_0\in L^2(\Omega)$ and let $u$ be the unique strong solution given by Theorem~$\ref{thmstrongH2}.$ Then,
\begin{gather}
\label{thmtdH21}
\vlim_{t\nearrow\infty}\|u(t)\|_{H^1_0(\Omega)}=\vlim_{t\nearrow\infty}\|u(t)\|_{L^p(\Omega)}=\vlim_{t\nearrow\infty}\frac{\d}{\d t}\|u(t)\|_{L^2(\Omega)}^2=0,
\end{gather}
for any $p\in\left(2m,\frac{2N}{N-2}\right]$ $(p\in(2m,\infty)$ if $N=2,$ $p\in(2m,\infty]$ if $N=1).$
\end{thm}

\begin{rmk}
\label{rmkthmtdH2}
Let the assumptions of Theorem~\ref{thmtdH2} be fulfilled. Below are some comments about the asymptotic behavior of the solution.
\begin{enumerate}
\item
\label{rmkthmtdH21}
If $m=0$ then $|\Omega|<\infty$ and by \eqref{thmtdH21}, $\vlim_{t\nearrow\infty}\|u(t)\|_{W^{1,q}(\Omega)}=0,$ for any $q\in(0,2].$
\item
\label{rmkthmtdH22}
Let $E=\big\{u\in H^1_0(\Omega); \Delta u\in L^2(\Omega)\big\}$ and $\|u\|_E^2=\|u\|_{L^2(\Omega)}^2+\|\Delta u\|_{L^2(\Omega)}^2,$ for any $u\in E.$ We recall that if $\Omega=\R^N,$ if $\Omega$ is a half-space or if $\Omega$ is bounded with a $C^{1,1}$-boundary then $E=H^2(\Omega)\cap H^1_0(\Omega)$ with equivalent norms. Indeed, this is due to Fourier's transform, Plancherel's formula, Haroske and Triebel~\cite[Theorem~5.16, p.131; Proposition~5.17, p.132]{MR2375667} and Grisvard~\cite[Corollary~2.2.2.4, p.91]{MR3396210}. With help of Property~\ref{thmstrongH24} of Theorem~\ref{thmstrongH2}, it follows from \eqref{thmtdH21} and Gagliardo-Nirenberg's inequality that,
\begin{gather*}
\vlim_{t\nearrow\infty}\|u(t)\|_{W^{1,q}(\Omega)}=\vlim_{t\nearrow\infty}\|u(t)\|_{L^p(\Omega)}=0,
\end{gather*}
for any $q\in\left[\frac43,2\right]$ with $q>\frac{4m}{m+1},$ and any $p\in\left(2m,\frac{2N}{N-4}\right)$ $(p\in(2m,\infty]$ if $N\le3).$
\end{enumerate}
\end{rmk}

\section{On the zero-order terms}
\label{functionals}

In this section we analize the functionals associated to the zero-order terms in equation~\eqref{nls}.

\begin{lem}
\label{lemVL}
Let $V=V_1+V_2\in L^\infty(\Omega;\R)+L^{p_V}(\Omega;\R),$ where $p_V$ is given by~\eqref{pV}. Then for any $u\in H^1_0(\Omega),$ we have $Vu\in L^2(\Omega)$ and,
\begin{gather}
\label{lemVL2}
\|Vu\|_{L^2(\Omega)}\le C\|V\|_{L^\infty(\Omega)+L^{p_V}(\Omega)}\|u\|_{H^1_0(\Omega)},
\end{gather}
where $C=C(N)$ $(C=C(\beta),$ if $N=2).$ In addition, for any $u\in H^1_0(\Omega),$
\begin{gather}
\label{lemV1}
\|V_1u\|_{L^2(\Omega)}\le\|V_1\|_{L^\infty(\Omega)}\|u\|_{L^2(\Omega)},
\end{gather}
and for any $\rho>0,$
\begin{gather}
\label{lemV2}
\|V_2u\|_{L^2(\Omega)}\le C\rho^{1-\gamma}\|V_2\|_{L^{p_V}(\Omega)}^{2-\gamma}\|u\|_{L^2(\Omega)}^\gamma+\frac1\rho\|\nabla u\|_{L^2(\Omega)}^2,
\end{gather}
where $\gamma=\frac23$ if $N=1,$ $\gamma=\frac\beta{\beta+1}$ if $N=2,$ $\gamma=0$ if $N\ge3,$ and $C=C(N)$ $\big(C=C(\beta),$ if $N=2\big).$
\end{lem}

\begin{lem}
\label{lemVH}
Let $V=V_1+V_2\in L^\infty(\Omega;\R)+L^{p_V}(\Omega;\R),$ where $p_V$ is given by~\eqref{pV}. Then for any $u\in L^2(\Omega),$ we get that $Vu\in H^{-1}(\Omega)$ and,
\begin{align}
\label{lemVH-1}
&	\|Vu\|_{H^{-1}(\Omega)}\le C\|V\|_{L^\infty(\Omega)+L^{p_V}(\Omega)}\|u\|_{L^2(\Omega)},	\\
\label{lemdualV}
&	\langle Vu,v\rangle_{H^{-1}(\Omega),H^1_0(\Omega)}=(u,Vv)_{L^2(\Omega)},
\end{align}
for any $v\in H^1_0(\Omega),$ where $C$ is given by~\eqref{lemVL2}.
\end{lem}

\begin{lem}
\label{lemcong}
The following properties are satisfied by the saturation terms $g_0^m(u):$
\begin{enumerate}
\item
\label{lemcong1}
Let $m\in(0,1].$ Then for any $p\in[1,\infty),$ we have that $g_0^m\in C\big(L^p(\Omega);L^\frac{p}m(\Omega)\big)$ and $g_0^m$ is bounded on bounded sets. More precisely,
\begin{gather*}
\|g_0^m(u)-g_0^m(v)\|_{L^\frac{p}m(\Omega)}\le3\|u-v\|_{L^p(\Omega)}^m,
\end{gather*}
for any $u,v\in L^p(\Omega).$
\item
\label{lemcong2}
Let $m\in[0,1]$ and $\eps>0.$ Then $g_\eps^m\in C\big(L^2(\Omega);L^2(\Omega)\big)$ and $g_\eps^m$ is bounded on bounded sets.
\end{enumerate}
\end{lem}

\begin{proof*}
The first part can be found in~Bégout and D\'{\i}az~\cite[Lemma~6.2]{MR4053613} while \ref{lemcong2} is obvious.
\medskip
\end{proof*}

\begin{vproof}{of Lemma~\ref{lemVL}.}
Let $u\in H^1_0(\Omega).$ By Hölder's inequality, we get~\eqref{lemV1} and,
\begin{gather}
\label{demlemV}
\|V_2u\|_{L^2(\Omega)}\le\|V_2\|_{L^{p_V}(\Omega)}\times
\begin{cases}
\|u\|_{L^\infty(\Omega)},					&	\text{if } N=1,	\medskip \\
\|u\|_{L^\frac{2(\beta+2)}\beta(\Omega)},		&	\text{if } N=2,	\smallskip \\
\|u\|_{L^\frac{2N}{N-2}(\Omega)},			&	\text{if } N\ge3.
\end{cases}
\end{gather}
Then \eqref{lemVL2} comes from \eqref{lemV1}, \eqref{demlemV} and the Sobolev embeddings. Let $\rho>0$ and $\nu=\rho\|V_2\|_{L^{p_V}(\Omega)}.$ By Sobolev's embedding and Gagliardo-Nirenberg's and Young's inequalities, we have
\begin{gather*}
\begin{cases}
\|u\|_{L^\infty(\Omega)}\le C\|u\|_{L^2(\Omega)}^\frac12\|\nabla u\|_{L^2(\Omega)}^\frac12
\le C\nu^\frac13\|u\|_{L^2(\Omega)}^\frac23+\frac1\nu\|\nabla u\|_{L^2(\Omega)}^2,								&	\text{if } N=1,	\medskip \\
\|u\|_{L^\frac{2(\beta+2)}\beta(\Omega)}\le C\|u\|_{L^2(\Omega)}^\frac\beta{\beta+2}\|\nabla u\|_{L^2(\Omega)}^\frac2{\beta+2}
\le C\nu^\frac1{\beta+1}\|u\|_{L^2(\Omega)}^\frac\beta{\beta+1}+\frac1\nu\|\nabla u\|_{L^2(\Omega)}^2,				&	\text{if } N=2,	\medskip \\
\|u\|_{L^\frac{2N}{N-2}(\Omega)}\le C\|\nabla u\|_{L^2(\Omega)}
\le C\nu+\frac1\nu\|\nabla u\|_{L^2(\Omega)}^2,															&	\text{if } N\ge3.
\end{cases}
\end{gather*}
Putting together \eqref{demlemV} and the above estimates, we obtain \eqref{lemV2}.
\medskip
\end{vproof}

\begin{vproof}{of Lemma~\ref{lemVH}.}
Let $u,v\in H^1_0(\Omega).$ By Lemma~\ref{lemVL}, $Vu\in L^2(\Omega)\inj H^{-1}(\Omega)$ with dense embedding and,
\begin{gather*}
\langle Vu,v\rangle_{H^{-1}(\Omega),H^1_0(\Omega)}=\langle Vu,v\rangle_{L^2(\Omega),L^2(\Omega)}
=(u,Vv)_{L^2(\Omega)},	\\
\sup_{\|v\|_{H^1_0(\Omega)}=1}\left|\langle Vu,v\rangle_{H^{-1}(\Omega),H^1_0(\Omega)}\right|
\le C\|V\|_{L^\infty(\Omega)+L^{p_V}(\Omega)}\|u\|_{L^2(\Omega)},
\end{gather*}
by Cauchy-Schwarz's inequality and \eqref{lemVL2}. The inequality is extended to any $u\in L^2(\Omega)$ by density. The lemma is proved.
\medskip
\end{vproof}

\section{Some maximal monotone operators}
\label{maxmon}

In all this section, we suppose Assumption~\ref{ass1} but with $a\in C(m),$ not merely $a\in C_{\mathrm{int}}(m),$ if $m\in(0,1)$ (unless if specified). Let $\eps\ge0.$
Let us define the following operators on $L^2(\Omega).$
\begin{align*}
&
\forall u\in D(L)\eqdef\left\{u\in H^1_0(\Omega); \; \Delta u\in L^2(\Omega)\right\}, \; Lu=-\vi\Delta u-\vi Vu,
\\[2pt]
&
\begin{cases}
D(B_\eps^m)=L^2(\Omega), \; \eps>0 \text{ or } m=1,	\medskip \\
\forall u\in D(B_\eps^m), \;  B_\eps^mu=-\vi ag_\eps^m(u),
\end{cases}
\begin{cases}
D(A_\eps^m)=D(L), \; \eps>0,	\medskip \\
\forall u\in D(A_\eps^m), \;  A_\eps^mu=Lu+B_\eps^mu,
\end{cases}
\end{align*}
\begin{align*}
&
\begin{cases}
D(B_0^0)=L^2(\Omega), \; |\Omega|<\infty,	\medskip \\
\forall u\in D(B_0^0), \;  B_0^0u=\Big\{U\in L^\infty(\Omega); \; \|U\|_{L^\infty(\Omega)}\le1 \text{ and if } u(x)\neq0, \; U(x)=g_0^0(u)(x)\Big\},
\end{cases}
\\[2pt]
&
\begin{cases}
D(A_0^0)=D(L), \; |\Omega|<\infty,	\medskip \\
\forall u\in D(A_0^0), \;  A_0^0u=\big\{Lu-\vi aU; \; U\in B_0^0u\big\},
\end{cases}
\\[2pt]
&
\begin{cases}
D(A_0^m)=\left\{u\in D(L); \; u\in L^{2m}(\Omega)\right\}, \; m>0,	\medskip \\
\forall u\in D(A_0^m), \;  A_0^mu=Lu-\vi ag_0^m(u).
\end{cases}
\end{align*}
It is clear that the all above domains are dense in $L^2(\Omega)$ since they all contain $\Dr(\Omega),$ which is dense in $L^2(\Omega).$

\begin{lem}
\label{lemenereg}
Let $u\in D(L)$ and $U\in B_0^0u.$ We have the following results.
\begin{enumerate}
\item
\label{lemenereg1}
If $m>0$ and if $u^m\Delta u\in L^1(\Omega)$ then $\Re\left(\vi a\dsp\int_\Omega g_0^m(u)\ovl{\Delta u}\d x\right)\ge0.$
\item
\label{lemenereg2}
If $m=0$ then $\Re\left(\vi a\dsp\int_\Omega U\ovl{\Delta u}\d x\right)\ge0.$
\end{enumerate}
\end{lem}

\begin{proof*}
By Bégout and D\'{\i}az~\cite[Lemma~6.3 and Remark~6.4]{MR4053613}, we only have to show \ref{lemenereg2}. Let $u\in D(L)$ and $U\in B_0^0u.$ Set $\omega=\big\{x\in\Omega;u(x)\neq0\big\}.$ Since $a\in\bigcap\limits_{0<m<1}C(m)$ and $g_0^m(u)\xrightarrow[m\searrow0]{\text{a.e.\,on }\omega}g_0^0(u),$ it follows fom the dominated convergence Theorem and \ref{lemenereg1} that,
\begin{gather}
\label{demlemenereg1}
\Re\left(\vi a\dsp\int_\omega\frac{u}{|u|}\ovl{\Delta u}\d x\right)\ge0.
\end{gather}
It is well-known that if $u\in H^1(\Omega)$ then $\nabla u=0,$ almost everywhere in $\omega^\co.$ In fact, since $u\in H^2_\loc(\Omega),$ $\Delta u=0,$ almost everywhere in $\omega^\co \cap K,$ for any compact subset $K\subset\Omega,$ hence in $\omega^\co.$ It follows that,
\begin{gather}
\label{demlemenereg2}
\Re\left(\vi a\dsp\int_{\omega^\co}U\ovl{\Delta u}\d x\right)=0.
\end{gather}
Summing \eqref{demlemenereg1} with \eqref{demlemenereg2}, we get the desired result.
\medskip
\end{proof*}

\begin{lem}
\label{lemL}
$(L,D(L))$ is a linear skew-adjoint operator on $L^2(\Omega)$ with dense domain. In particular, it is maximal monotone.
\end{lem}

\begin{proof*}
It is clear that $Lu\in L^2(\Omega),$ for any $u\in D(L)$ (Lemma~\ref{lemVL}) and that $(L,D(L))$ is a skew-adjoint linear operator with dense domain, from which the result follows (Cazenave and Haraux~\cite[Corollary~2.4.9, p.24]{MR2000e:35003}).
\medskip
\end{proof*}

\noindent
The monotonicity result below is a slight generalization of a result of Hayashi~\cite[Lemma~4.3]{MR3802567} but for the convenience of the reader, we give its proof. Actually, in his paper the quantity in~\eqref{lemmongen1} below is nonegative and we need a positive quantity.

\begin{lem}
\label{lemmongen}
Let $f:(0,\infty)\tends(0,\infty)$ be an increasing function. Then for any $(z_1,z_2)\in\C^2$ such that $z_1z_2\neq0$ and $|z_1|\neq|z_2|,$
\begin{gather}
\label{lemmongen1}
\Re\left(\left(f(|z_1|)\frac{z_1}{|z_1|}-f(|z_2|)\frac{z_2}{|z_2|}\right)\big(\ovl{z_1-z_2}\big)\right)>0.
\end{gather}
If $f$ is merely nondecreasing or if $|z_1|=|z_2|$ then the quantity in~\eqref{lemmongen1} is nonnegative.
\end{lem}

\begin{proof*}
Let $f:(0,\infty)\tends(0,\infty)$ be an increasing function and let $(z_1,z_2)\in\C^2$ be such that $z_1z_2\neq0$ and $|z_1|\neq|z_2|.$ We have,
\begin{align*}
	&	\; \Re\left(\left(f(|z_1|)\frac{z_1}{|z_1|}-f(|z_2|)\frac{z_2}{|z_2|}\right)\big(\ovl{z_1-z_2}\big)\right)			\\
   =	&	\; f(|z_1|)|z_1|-f(|z_1|)\frac{\Re(z_1\ovl{z_2})}{|z_1|}-f(|z_2|)\frac{\Re(z_1\ovl{z_2})}{|z_2|}+f(|z_2|)|z_2|	\\
 \ge	&	\; f(|z_1|)|z_1|-f(|z_1|)|z_2|-f(|z_2|)|z_1|+f(|z_2|)|z_2|											\\
  =	&	\; \big(f(|z_1|-f(|z_2|)\big)\big(|z_1|-|z_2|\big)												\\
  >	&	\; 0,
\end{align*}
since $f$ is increasing and $|z_1|\neq|z_2|.$
\medskip
\end{proof*}

\begin{rmk}
\label{rmklemmongen}
Since on $\C\setminus\{0\},$ $\Re(z_1\ovl{z_2})=|z_1||z_2|$ if, and only if, $\Arg(z_1)=\Arg(z_2),$ it follows from the proof of Lemma~\ref{lemmongen} that if $f:(0,\infty)\tends(0,\infty)$ is a nondecreasing function then for any $(z_1,z_2)\in\C^2$ such that $z_1z_2\neq0$ and $\Arg(z_1)\neq\Arg(z_2)$ (with possibly $|z_1|=|z_2|),$ the quantity in~\eqref{lemmongen1} is positive (and not merely nonnegative). Here and after, $\Arg(z)\in(-\pi,\pi]$ denotes the principal value of the argument of $z\in\C\setminus\{0\}.$

\end{rmk}

\begin{cor}
\label{cormongen}
Let $(z_1,z_2)\in\C^2.$
\begin{enumerate}
\item
\label{cormongen1}
Assume that $m+\eps>0.$ If $|z_1|\neq|z_2|$ then,
\begin{gather}
\label{cormongen11}
\Re\left(\left((|z_1|^2+\eps)^\frac{m-1}2z_1-(|z_2|^2+\eps)^\frac{m-1}2z_2\right)\big(\ovl{z_1-z_2}\big)\right)>0,
\end{gather}
$($respectively, $\ge0,$ if $|z_1|=|z_2|).$
\item
\label{cormongen2}
Assume that $m+\eps=0.$ If $z_1z_2\neq0$ then,
\begin{gather}
\label{cormongen21}
\Re\left(\left(\frac{z_1}{|z_1|}-\frac{z_2}{|z_2|}\right)\big(\ovl{z_1-z_2}\big)\right)\ge0.
\end{gather}
\end{enumerate}
\end{cor}

\begin{proof*}
Apply Lemma~\ref{lemmongen}, where for any $t>0,$ $f(t)=(t^2+\eps)^\frac{m-1}2t.$
\medskip
\end{proof*}

\noindent
The result below, for $\eps=0,$ is due to Liskevich and Perel$^\p$muter~\cite[Lemma~2.2]{MR1224619}. Nevertheless, we will need to generalize it to the regularized case $\eps>0.$

\begin{lem}
\label{lemmon}
We have,
\begin{align*}
	&	\; 2\sqrt m\left|\Im\left(\left((|z_1|^2+\eps)^\frac{m-1}2z_1-(|z_2|^2+\eps)^\frac{m-1}2z_2\right)\big(\ovl{z_1-z_2}\big)\right)\right|	\\
  \le	&	\; (1-m)\Re\left(\left((|z_1|^2+\eps)^\frac{m-1}2z_1-(|z_2|^2+\eps)^\frac{m-1}2z_2\right)\big(\ovl{z_1-z_2}\big)\right).
\end{align*}
for any $(z_1,z_2)\in\C\times\C$ $($and $z_1z_2\neq0,$ if $m=\eps=0).$
\end{lem}

\begin{rmk}
\label{rmklemmon}
If $m=0$ then Lemma~\ref{lemmon} is nothing else but Corollary~\ref{cormongen} (while if $m=1$ then the conclusion is that the complex number we are computing between the parentheses is a nonnegative real number, which is obvious).
\end{rmk}

\begin{vproof}{of Lemma~\ref{lemmon}.}
By Remark~\ref{rmklemmon}, we may assume that $0<m<1.$ Let $(z_1,z_2)\in\C^2.$ \\
Set $Z_\eps=\left((|z_1|^2+\eps)^\frac{m-1}2z_1-(|z_2|^2+\eps)^\frac{m-1}2z_2\right)\big(\ovl{z_1-z_2}\big).$ A straightforward calculation gives,
\begin{align*}
&	\Re(Z_\eps)=\left(|z_1|^2(|z_1|^2+\eps)^\frac{m-1}2+|z_2|^2(|z_2|^2+\eps)^\frac{m-1}2\right)
		-\Re(z_1\ovl{z_2})\left(|z_1|^2+\eps)^\frac{m-1}2+(|z_2|^2+\eps)^\frac{m-1}2\right),										\\
&	\Im(Z_\eps)=\Im(\ovl{z_1}z_2)\left((|z_1|^2+\eps)^\frac{m-1}2-(|z_2|^2+\eps)^\frac{m-1}2\right),									\\
&	\Re(z_1\ovl{z_2})=|z_1|\,|z_2|\cos\big(\Arg(z_1\ovl{z_2})\big),	\quad	\Im(\ovl{z_1}z_2)=|z_1|\,|z_2|\sin\big(\Arg(\ovl{z_1}z_2)\big).
\end{align*}
Note that $\Im(Z_\eps)=0\le\Re(Z_\eps)$ if $z_1z_2=0$ or $|z_1|=|z_2|$ (Corollary~\ref{cormongen}). So we may assume that $|z_1|>|z_2|>0.$ We set $t=|z_1|,$ $s=|z_2|$ and $\theta=\Arg(\ovl{z_1}z_2).$ By Corollary~\ref{cormongen}, $\Re(Z_\eps)>0$ and we may define $F_\eps$ by,
\begin{gather*}
F_\eps(t,s,\theta)=\frac{|\Im(Z_\eps)|}{\Re(Z_\eps)}.
\end{gather*}
Since $F_\eps\ge0,$ we shall show that,
\begin{gather}
\label{demlemmon1}
F_\eps(t,s,\theta)^2\le\frac{(1-m)^2}{4m},
\end{gather}
with,
\begin{align*}
F_\eps(t,s,\theta)^2
	=	& \; \frac{t^2s^2\left((t^2+\eps)^\frac{m-1}2-(s^2+\eps)^\frac{m-1}2\right)^2(1-\cos^2\theta)}
			{\left(\left(t^2(t^2+\eps)^\frac{m-1}2+s^2(s^2+\eps)^\frac{m-1}2\right)
			-ts\left((t^2+\eps)^\frac{m-1}2+(s^2+\eps)^\frac{m-1}2\right)\cos\theta\right)^2},			\\
\eqdef	& \; \frac{A(1-\cos^2\theta)}{(B-C\cos\theta)^2}.
\end{align*}
We proceed with the proof in four steps. 
\\
\textbf{Step 1:} $F_\eps(t,s,\theta)^2\le\dfrac{t^2s^2\left((t^2+\eps)^\frac{m-1}2-(s^2+\eps)^\frac{m-1}2\right)^2}
{(t^2-s^2)\big(t^2(t^2+\eps)^{m-1}-s^2(s^2+\eps)^{m-1}\big)}.$ \\
We write $\sigma=\cos\theta$ and $g(\sigma)=F_\eps(t,s,\theta)^2=\frac{A(1-\sigma^2)}{(B-C\sigma)^2}.$ Note that since $t>s>0$ then, with help of Corollary~\ref{cormongen}, we have $A>0,$ $B>0$ and $B-C\sigma>0,$ for any $\sigma\in[-1,1].$ In particular, $0<C<B.$ A study of $g$ gives,
\begin{gather*}
\max_{\sigma\in[-1,1]}g(\sigma)=g\left(\frac{C}{B}\right).
\end{gather*}
It follows that, $\vsup_{\theta\in(-\pi,\pi]}F_\eps(t,s,\theta)^2\le g\left(\frac{C}{B}\right),$ which gives the desired result. \\
\textbf{Step 2:} $\left((t^2+\eps)^\frac{m-1}2-(s^2+\eps)^\frac{m-1}2\right)^2\le\dfrac{(1-m)^2}{4m}\dfrac{t^2-s^2}{(t^2+\eps)(s^2+\eps)}\big((t^2+\eps)^m-(s^2+\eps)^m\big).$ \\
By the Cauchy-Schwarz inequality, we have,
\begin{align*}
	&	\; \left((t^2+\eps)^\frac{m-1}2-(s^2+\eps)^\frac{m-1}2\right)^2											\\
   =	&	\; \frac{(1-m)^2}4\left(\:\vint_{s^2+\eps}^{t^2+\eps}\sigma^\frac{m-3}2\d\sigma\right)^2
			=\frac{(1-m)^2}4\left(\:\vint_{s^2+\eps}^{t^2+\eps}\sigma^{-1}\sigma^\frac{m-1}2\d\sigma\right)^2			\\
  \le	&	\; \frac{(1-m)^2}4\vint_{s^2+\eps}^{t^2+\eps}\sigma^{-2}\d\sigma\vint_{s^2+\eps}^{t^2+\eps}\sigma^{m-1}\d\sigma	\\
   =	&	\; \frac{(1-m)^2}{4m}\big((s^2+\eps)^{-1}-(t^2+\eps)^{-1}\big)\big((t^2+\eps)^m-(s^2+\eps)^m\big),
\end{align*}
which is Step~2.
\\
\textbf{Step 3:} $0<(t^2+\eps)^m-(s^2+\eps)^m\le t^2(t^2+\eps)^{m-1}-s^2(s^2+\eps)^{m-1}.$ \\
Indeed,
\begin{align*}
	&	\; (t^2+\eps)^m-(s^2+\eps)^m	\\
   =	&	\; \left(t^2(t^2+\eps)^{m-1}-s^2(s^2+\eps)^{m-1}\right)-\eps\left((s^2+\eps)^{m-1}-(t^2+\eps)^{m-1}\right)	\\
  \le	&	\; t^2(t^2+\eps)^{m-1}-s^2(s^2+\eps)^{m-1},
\end{align*}
since $t>s>0$ and $m-1<0.$ Hence Step~3.
\\
\textbf{Step 4:} Conclusion. \\
Putting together Steps~1--3, we infer,
\begin{gather*}
F_\eps(t,s,\theta)^2\le\frac{(1-m)^2}{4m}\frac{t^2s^2}{(t^2+\eps)(s^2+\eps)}\le\frac{(1-m)^2}{4m},
\end{gather*}
which is \eqref{demlemmon1}. This ends the proof.
\medskip
\end{vproof}

\begin{cor}
\label{corlemineg}
Assume $m+\eps>0.$ Let $u,v\in L^{m+1}(\Omega)$ if $\eps=0,$ and let $u,v\in L^2(\Omega)$ if $\eps>0.$ Then $\big(g_\eps^m(u)-g_\eps^m(v)\big)(\ovl{u-v})\in L^1(\Omega)$ and,
\begin{gather}
\label{corlemineg1}
\Re\left(-\vi a\vint_\Omega\big(g_\eps^m(u)-g_\eps^m(v)\big)(\ovl{u-v})\d x\right)\ge0,
\end{gather}
for any $a\in C(m).$
\end{cor}

\begin{proof*}
Assume $m\in[0,1)$ and $\eps\ge0$ with $m+\eps>0.$ Let $u,v$ be as in the corollary. Then by Lemma~\ref{lemcong} and Hölder's inequality, $\big(g_\eps^m(u)-g_\eps^m(v)\big)(\ovl{u-v})\in L^1(\Omega).$ Now, let $a\in C(m).$ By Lemma~\ref{lemmon},
\begin{align*}
	& \; \Re\left(-\vi a\vint_\Omega\big(g_\eps^m(u)-g_\eps^m(v)\big)(\ovl{u-v})\d x\right)								\\
   =	& \; \Im(a)\Re\vint_\Omega\big(g_\eps^m(u)-g_\eps^m(v)\big)\big(\ovl{u-v}\big)\d x
   		+\Re(a)\Im\vint_\Omega\big(g_\eps^m(u)-g_\eps^m(v)\big)\big(\ovl{u-v}\big)\d x								\\
 \ge	& \; \left(\Im(a)-|\Re(a)|\frac{1-m}{2\sqrt m}\right)\Re\vint_\Omega\big(g_\eps^m(u)-g_\eps^m(v)\big)\big(\ovl{u-v}\big)\d x	\\
 \ge	& \; 0,
\end{align*}
if $m>0.$ If $m=0$ then $a\in C(0)=\{0\}\times\vi(0,\infty)$ and,
\begin{gather*}
\Re\left(-\vi a\vint_\Omega\big(g_\eps^0(u)-g_\eps^0(v)\big)(\ovl{u-v})\d x\right)
=\Im(a)\Re\vint_\Omega\big(g_\eps^0(u)-g_\eps^0(v)\big)\big(\ovl{u-v}\big)\d x\ge0,
\end{gather*}
by Corollary~\ref{cormongen}. This ends the proof.
\medskip
\end{proof*}

\begin{cor}
\label{corAmon}
Assume $m\in(0,1]$ and $a\in C(m).$ Then $(A_0^m,D(A_0^m))$ is monotone on $L^2(\Omega)$ with dense domain.
\end{cor}

\begin{proof*}
By Lemmas~\ref{lemL} and \ref{lemcong}, $A_0^m:D(A_0^m)\tends L^2(\Omega)$ is well-defined. Let $u,v\in D(A_0^m).$ We have, $D(A_0^m)\subset L^{2m}(\Omega)\cap L^2(\Omega)\subset L^{m+1}(\Omega),$ and so Corollary~\ref{corlemineg} applies. Finally, by skew-adjointness of $L$ (Lemma~\ref{lemL}),
\begin{gather*}
(A_0^mu-A_0^mv,u-v)_{L^2(\Omega)}=\Re\left(-\vi a\vint_\Omega\big(g_0^m(u)-g_0^m(v)\big)(\ovl{u-v})\d x\right)\ge0,
\end{gather*}
by Corollary~\ref{corlemineg}.
\medskip
\end{proof*}

\begin{cor}
\label{corA0mon}
Assume $a\in C(0).$ Then $(A_0^0,D(A_0^0))$ is monotone on $L^2(\Omega)$ with dense domain.
\end{cor}

\begin{proof*}
Since $|\Omega|<\infty,$ $L^\infty(\Omega)\inj L^2(\Omega)$ and we have $A_0^0u\in\vPr\big(L^2(\Omega)\big),$ for any $u\in D(A_0^0).$ Since $a\in C(0),$ we have $a=\vi\lambda,$ for some real $\lambda>0.$ Let $u_1,u_2\in D(A_0^0)$ and $(V_1,V_2)\in A_0^0u_1\times A_0^0u_2.$ Then for each $j\in\{1,2\},$ there exists $U_j\in B_0^0u_j$ such that $V_j=L u_j+\lambda U_j.$ By skew-adjointness of $L,$
\begin{gather*}
(V_1-V_2,u_1-u_2)_{L^2(\Omega)}=\lambda(U_1-U_2,u_1-u_2)_{L^2(\Omega)}.
\end{gather*}
For each $j\in\{1,2\},$ we define, $\omega_j=\big\{x\in\Omega; \; u_j(x)\neq0\big\}.$ We then have,
\begin{align*}
	&	\; (U_1-U_2,u_1-u_2)_{L^2(\Omega)}													\\
   =	&	\; \Re\left(\:\vint_{\omega_1^\co\cap\,\omega_2}\left(U_1-\frac{u_2}{|u_2|}\right)\ovl{(-u_2)}\d x\right)
			+\Re\left(\:\vint_{\omega_1\cap\,\omega_2^\co}\left(\frac{u_1}{|u_1|}-U_2\right)\ovl{u_1}\d x\right)	\\
	&	\qquad\qquad \; +\Re\left(\:\vint_{\omega_1\cap\,\omega_2}
			\left(\frac{u_1}{|u_1|}-\frac{u_2}{|u_2|}\right)\big(\ovl{u_1-u_2}\big)\d x\right)					\\
  \ge	&	\; \Re\left(\:\vint_{\omega_1^\co\cap\,\omega_2}\big(|u_2|-U_1\ovl{u_2}\big)\d x\right)
			+\Re\left(\:\vint_{\omega_1\cap\,\omega_2^\co}\big(|u_1|-U_2\ovl{u_1}\big)\d x\right) 				\\
  \ge &	\; 0.
\end{align*}
Indeed, the first inequality is due to~\eqref{cormongen21}, while the last one comes from the fact that $|U_1\ovl{u_2}|\le|u_2|$ and $|U_2\ovl{u_1}|\le|u_1|.$ This ends the proof.
\medskip
\end{proof*}

\begin{cor}
\label{corAe}
Assume $m\in[0,1)$ and $\eps>0,$ or $(m,\eps)=(1,0).$ Let $a\in C(m).$ Then $(A_\eps^m,D(A_\eps^m))$ is maximal monotone on $L^2(\Omega)$ with dense domain.
\end{cor}

\begin{proof*}
By Lemma~\ref{lemL}, $(L,D(L))$ is maximal monotone and by Lemma~\ref{lemcong}, $D(B_\eps^m)=L^2(\Omega),$ $B_\eps^m\in C\big(L^2(\Omega);L^2(\Omega)\big)$ and
\begin{gather*}
(B_\eps^mu-B_\eps^mv,u-v)_{L^2(\Omega)}=\Re\left(-\vi a\vint_\Omega\big(g_\eps^m(u)-g_\eps^m(v)\big)(\ovl{u-v})\d x\right)\ge0,
\end{gather*}
for any $u,v\in L^2(\Omega)$ (Corollary~\ref{corlemineg}). We then deduce that $(B_\eps^m,L^2(\Omega))$ is maximal monotone (Brezis~\cite[Corollary~2.5, p.33]{MR0348562}) and so is, from abstract perturbations results, $A_\eps^m\eqdef L+B_\eps^m$ (Brezis~\cite[Corollary~2.7, p.36]{MR0348562}).
\medskip
\end{proof*}

\begin{lem}
\label{lemAmax}
Assume $m\in(0,1)$ and $a\in C_{\mathrm{int}}(m),$ or $m=0$ and $a\in C(0).$ Then, $R(I+A_0^m)=L^2(\Omega).$
\end{lem}

\begin{proof*}
Let $F\in L^2(\Omega).$ We proceed with the proof in five steps. \\
\textbf{Step 1:} Let $\eps>0.$ There exists $u_\eps\in D(A_\eps^m)$ satisfying,
\begin{gather}
\label{lemAmax1}
-\vi\Delta u_\eps-\vi Vu_\eps-\vi ag_\eps^m(u_\eps)+u_\eps=F, \text{ in } L^2(\Omega).
\end{gather}
Since $(A_\eps^m,D(A_\eps^m))$ is maximal monotone (Corollary~\ref{corAe}), we have $R(I+A_\eps^m)=L^2(\Omega)$ (Brezis~\cite[Proposition~2.2, p.23]{MR0348562}).
\\
\textbf{Step 2:} The families $(u_\eps)_{\eps>0}$ and $(Vu_\eps)_{\eps>0}$ are bounded in $H^1_0(\Omega)$ and in $L^2(\Omega),$ respectively, and there exist a $u\in H^1_0(\Omega)$ and a decreasing sequence $(\eps_n)_{n\in\N}\subset(0,\infty)$ converging toward $0$ such that $Vu\in L^2(\Omega)$ and,
\begin{align}
\label{lemAmax2}
&	u_{\eps_n}\xrightarrow[n\to\infty]{\Dr^\p(\Omega)}u,		\\
\label{lemAmax3}
&	Vu_{\eps_n}\xrightarrow[n\to\infty]{\Dr^\p(\Omega)}Vu,	\\
\label{lemAmax4}
&	u_{\eps_n}\xrightarrow[n\to\infty]{L^2_\loc(\Omega)}u,	\\
\label{lemAmax5}
&	u_{\eps_n}\xrightarrow[n\to\infty]{\text{a.e.\,in }\Omega}u.
\end{align}
Let $\eps>0.$ We successively take the $L^2$-scalar product of~\eqref{lemAmax1} with $u_\eps$ and then with $\vi u_\eps.$ We get,
\begin{gather}
\label{lemAmax6}
\Im(a)\vint_\Omega(|u_\eps|^2+\eps)^{-\frac{1-m}2}|u_\eps|^2\d x+\|u_\eps\|_{L^2(\Omega)}^2=\Re\vint_\Omega F\,\ovl{u_\eps}\d x,		\\
\label{lemAmax7}
\|\nabla u_\eps\|_{L^2(\Omega)}^2-\vint_\Omega V|u_\eps|^2\d x-\Re(a)\vint_\Omega(|u_\eps|^2+\eps)^{-\frac{1-m}2}|u_\eps|^2\d x
=\Im\vint_\Omega F\,\ovl{u_\eps}\d x.
\end{gather}
Applying Cauchy-Schwarz's inequality to~\eqref{lemAmax6}, we obtain $\|u_\eps\|_{L^2(\Omega)}\le\|F\|_{L^2(\Omega)}$ and so,
\begin{gather}
\label{lemAmax8}
\Im(a)\vint_\Omega(|u_\eps|^2+\eps)^{-\frac{1-m}2}|u_\eps|^2\d x+\|u_\eps\|_{L^2(\Omega)}^2\le\|F\|_{L^2(\Omega)}^2,
\end{gather}
Using Hölder's and Cauchy-Schwarz's inequalities in~\eqref{lemAmax7}, we get by \eqref{lemAmax8},
\begin{gather}
\label{lemAmax9}
\|\nabla u_\eps\|_{L^2(\Omega)}^2\le\left(1+\frac{\Re(a)_+}{\Im(a)}\right)\|F\|_{L^2(\Omega)}^2+\|Vu_\eps\|_{L^2(\Omega)}\|F\|_{L^2(\Omega)}.
\end{gather}
Let us write $V=V_1+V_2$ with $(V_1,V_2)\in L^\infty(\Omega;\R)\times L^{p_V}(\Omega;\R).$ Then by~\eqref{lemV1}, \eqref{lemV2} and \eqref{lemAmax8},
\begin{gather}
\label{lemAmax10}
\|Vu_\eps\|_{L^2(\Omega)}\le C\left(\|V_1\|_{L^\infty(\Omega)}+\|V_2\|_{L^{p_V}(\Omega)}^{2-\gamma}\right)\|F\|_{L^2(\Omega)}+\frac1{2\|F\|_{L^2(\Omega)}}\|\nabla u_\eps\|_{L^2(\Omega)}^2,
\end{gather}
where $C=C(N)$ $(C=C(\beta),$ if $N=2).$ Putting together~\eqref{lemAmax9} and \eqref{lemAmax10}, we infer
\begin{gather}
\label{lemAmax11}
\sup_{\eps>0}\|\nabla u_\eps\|_{L^2(\Omega)}+\sup_{\eps>0}\|Vu_\eps\|_{L^2(\Omega)}<\infty.
\end{gather}
By \eqref{lemAmax8} and \eqref{lemAmax11}, $(u_\eps)_{\eps>0}$ and $(Vu_\eps)_{\eps>0}$ are bounded in $H^1_0(\Omega)$ and in $L^2(\Omega),$ respectively. Since both spaces are reflexive, we obtain \eqref{lemAmax2}--\eqref{lemAmax5} for some $u\in H^1_0(\Omega)$ with $Vu\in L^2(\Omega)$ by local compactness, \eqref{lemVL2}, \eqref{lemdualV} and a decreasing sequence $\eps_n\searrow0.$
\\
\textbf{Step 3:} $u\in D(A_0^m)$ and if $m=0$ then $\vsup_{n\in\N}\|g_{\eps_n}^0(u_{\eps_n})\|_{L^\infty(\Omega)}\le1.$
\\
If $m=1$ then the result is a direct consequence of Step~2, \eqref{lemVL2} and the equation~\eqref{lemAmax1}. We continue with the case $m>0.$ Since $a\in C_{\mathrm{int}}(m),$ there exists $b\in\C$ such that $|b|=1,$ $\Re(b)>0,$ $\Im(b)<0$ and $ab\in C_{\mathrm{int}}(m)$ (Bégout~\cite[Lemma~4.2]{MR4098330}). We take the $L^2$-scalar product of~\eqref{lemAmax1} with $abg_\eps^m(u_\eps).$ We then get,
\begin{multline*}
\Re\left(\vi ab\vint_\Omega g_\eps^m(u_\eps)\ovl{\Delta u_\eps}\d x\right)-\Im(ab)\vint_\Omega Vg_\eps^m(u_\eps)\ovl{u_\eps}\d x
														+|a|^2\,|\Im(b)|\,\|g_\eps^m(u_\eps)\|_{L^2(\Omega)}^2	\\
+\Re\left(ab\vint_\Omega g_\eps^m(u_\eps)\ovl{u_\eps}\d x\right)=\Re\left(ab\vint_\Omega\ovl F\,g_\eps^m(u_\eps)\d x\right),
\end{multline*}
By~(6.8) in~Bégout and D\'{\i}az~\cite[Lemma~6.3]{MR4053613}, the first term in the left hand side of the above equality is nonnegative. With help of Cauchy-Schwarz's and Young's inequalities, and Step~2, we infer,
\begin{align*}
	&	\;	|a|\,|\Im(b)|\,\|g_\eps^m(u_\eps)\|_{L^2(\Omega)}^2								\\
  \le	&	\; \left(\sup_{\eps>0}\|Vu_\eps\|_{L^2(\Omega)}+\sup_{\eps>0}\|u_\eps\|_{L^2(\Omega)}
							+\|F\|_{L^2(\Omega)}\right)\|g_\eps^m(u_\eps)\|_{L^2(\Omega)}		\\
  \le	&	\; C+\frac{|a|\,|\Im(b)|}2\|g_\eps^m(u_\eps)\|_{L^2(\Omega)}^2,
\end{align*}
and thus $\sup_{\eps>0}\|g_\eps^m(u_\eps)\|_{L^2(\Omega)}<\infty.$ By~\eqref{lemAmax1} and Step~2, we deduce that,
\begin{gather*}
\sup_{\eps>0}\|\Delta u_\eps\|_{L^2(\Omega)}+\sup_{\eps>0}\|g_\eps^m(u_\eps)\|_{L^2(\Omega)}<\infty.
\end{gather*}
This last estimate with Step~2 and Fatou's Lemma imply that $\Delta u\in L^2(\Omega)$ and $g_0^m(u)\in L^2(\Omega).$ This last point means that $u\in L^{2m}(\Omega)$ and finally $u\in D(A_0^m).$ Now, we turn out to the case $m=0.$ In particular, $|\Omega|<\infty.$ We have $g_\eps^0(u_\eps)(x)=0,$ if $u_\eps(x)=0$ and $|g_\eps^0(u_\eps)(x)|\le|g_0^0(u_\eps)(x)|=1,$ otherwise. With the embedding $L^\infty(\Omega)\inj L^2(\Omega),$ Step~2 and \eqref{lemAmax1}, this implies that $(\Delta u_\eps)_{\eps>0}$ is bounded in $L^2(\Omega).$ Hence $u\in D(A_0^0)$ by \eqref{lemAmax2}.
\\
\textbf{Step 4:} If $m=0$ then there exists $U\in B_0^0u$ such that, up to a subsequence, $g_{\eps_n}^0(u_{\eps_n})\xrightarrow[n\to\infty]{\Dr^\p(\Omega)}U.$
\\
For any $n\in\N,$ $\|g_{\eps_n}^0(u_{\eps_n})\|_{L^\infty(\Omega)}\le1$ and by \eqref{lemAmax5}, $g_{\eps_n}^0(u_{\eps_n})(x)\xrightarrow[n\to\infty]{}g_0^0(u)(x),$ for almost every $x\in\Omega$ such that $u(x)\neq0.$ Then applying~Cazenave~\cite[Proposition~1.2.1, p.3]{MR2002047}, we get the desired result.
\\
\textbf{Step~5:} Conclusion.
\\
By~\eqref{lemAmax1} and Steps~2--4, if $m=0$ then for some $U\in B_0^0u,$ $u-\vi\Delta u-\vi Vu-\vi aU=F,$ in $\Dr^\p(\Omega),$ so in $L^2(\Omega),$ since $u\in D(A_0^0).$ In other words,
\begin{gather*}
u\in D(A_0^0)	\; \text{ and } \;	(I+A_0^0)u\ni F.
\end{gather*}
This ends the proof for $m=0.$ Now, assume that $m>0.$ By Step~3, $u\in D(A_0^m).$ It remains to show that, $(I+A_0^m)u=F.$ Let $\vphi\in\Dr(\Omega).$ By~\eqref{lemAmax1}, we have for any $n\in\N,$
\begin{gather}
\label{lemAmax12}
\langle u_{\eps_n}-\vi\Delta u_{\eps_n}-\vi Vu_{\eps_n},\vphi\rangle_{\Dr^\p(\Omega),\Dr(\Omega)}
-\Re\left(\vi a\vint_\Omega g_{\eps_n}^m(u_{\eps_n})\ovl\vphi\,\d x\right)=\langle F,\vphi\rangle_{\Dr^\p(\Omega),\Dr(\Omega)}.
\end{gather}
Let $\Omega^\p$ a bounded open subset of $\R^N$ be such that $\supp\vphi\subset\Omega^\p\subset\Omega.$ By~\eqref{lemAmax4}, there exist $h\in L^2(\Omega^\p;\R)$ and a subsequence, that we still denote by $(\eps_n)_{n\in\N},$ such that for any $n\in\N,$ $|u_{\eps_n}|\le h,$ almost everywhere in $\Omega^\p$ (see, for instance, Brezis~\cite[Theorem~4.9, p.94]{MR2759829}). Extending $h$ by $0$ over $\Omega\setminus\Omega^\p$ (with no change of notation), we obtain that $h\in L^2(\Omega;\R).$ With help of~\eqref{lemAmax5}, we obtain,
\begin{gather*}
g_{\eps_n}^m(u_{\eps_n})\,\ovl\vphi\xrightarrow[n\to\infty]{\text{a.e. in }\Omega}g_0^m(u)\,\ovl\vphi,			\\
|g_{\eps_n}^m(u_{\eps_n})\,\ovl\vphi|\le h^m|\vphi|, \text{ a.e. in } \Omega,
\end{gather*}
for any $n\in\N.$ But $h^m|\vphi|\in L^1(\Omega;\R)$ by Hölder's inequality. Applying the dominated convergence Theorem, we may pass to the limit in~\eqref{lemAmax12} to get with help of~\eqref{lemAmax2} and~\eqref{lemAmax3},
\begin{gather*}
u-\vi\Delta u-\vi Vu-\vi ag_0^m(u)=F, \text{ in } \Dr^\p(\Omega).
\end{gather*}
But $u\in D(A_0^m)$ and so the above equation makes sense in $L^2(\Omega).$ We conclude that,
\begin{gather*}
u\in D(A_0^m)	\; \text{ and } \;	(I+A_0^m)u=F, \text{ in } L^2(\Omega).
\end{gather*}
This ends the proof.
\medskip
\end{proof*}

\begin{cor}
\label{corAmaxmon}
Assume $m\in(0,1)$ and $a\in C_{\mathrm{int}}(m),$ or $m\in\{0,1\}$ and $a\in C(m).$ Then $(A_0^m,D(A_0^m))$ is maximal monotone on $L^2(\Omega)$ with dense domain.
\end{cor}

\begin{proof*}
If $m=1$ then the result comes from Corollary~\ref{corAe}. Now assume that $0\le m<1.$ Since $(A_0^m,D(A_0^m))$ is monotone (with dense domain) and $R(I+A_0^m)=L^2(\Omega)$ (Corollary~\ref{corAmon}, Corollary~\ref{corA0mon} and Lemma~\ref{lemAmax}), $(A_0^m,D(A_0^m))$ is maximal monotone (Brezis~\cite[Proposition~2.2, p.23]{MR0348562}).
\medskip
\end{proof*}

\section{Proofs of the existence theorems}
\label{proofexi}

In this section, we shall use the notations of the previous section.
\medskip

\begin{vproof}{of Proposition~\ref{propsolL2}.}
Set $Y=H^2_0(\Omega)\cap L^\frac2{2-m}(\Omega).$ Then, $Y^\star=H^{-2}(\Omega)+L^\frac2m(\Omega).$
By \eqref{cv}, \eqref{lemVH-1} and Lemma~\ref{lemcong},
\begin{align}
\label{dempropsolL21}
&	\Delta u_n\xrightarrow[n\to\infty]{C([0,T];H^{-2}(\Omega))}\Delta u,		\\
\label{dempropsolL22}
&	Vu_n\xrightarrow[n\to\infty]{C([0,T];H^{-1}(\Omega))}Vu,				\\
\label{dempropsolL23}
&	g(u_n)\xrightarrow[n\to\infty]{C([0,T];L^\frac2m(\Omega))}g(u), \text{ if } m>0,
\end{align}
for any $T>0.$ If $m=0$ then by Definition~\ref{defsol}, $\vsup_{n\in\N}\|U_n\|_{L^\infty((0,\infty)\times\Omega)}\le1.$ By \eqref{cv}, up to a subsequence, $u_n\xrightarrow[n\to\infty]{\text{a.e.\,in }(0,\infty)\times\Omega}u.$ By the Vitali Theorem, there exist a subsequence $(U_{n_k})_{k\in\N}\subset(U_n)_{n\in\N}$ and $U\in L^\infty\big((0,\infty)\times\Omega\big)$ such that,
\begin{align}
\label{dempropsolL24}
&	U_{n_k}\underset{k\to\infty}{\overset{L^\infty_{\w\star}((0,\infty)\times\Omega)}{-\!\!\!-\!\!\!-\!\!\!-\!\!\!-\!\!\!-\!\!\!-\!\!\!-\!\!\!-\!\!\!-\!\!\!-\!\!\!\weak}}U,	\\
\label{dempropsolL25}
&	U(t,x)=\dfrac{u(t,x)}{|u(t,x)|}, \text{ if } u(t,x)\neq0,																		\\
\label{dempropsolL26}
&	\|U\|_{L^\infty((0,\infty)\times\Omega)}\le1.
\end{align}
Then it follows from the equation satisfied by $u_n,$ \eqref{cv} and \eqref{dempropsolL21}--\eqref{dempropsolL26} that \eqref{propsolL21} holds true and $u$ solves~\eqref{nls} in $L^1_\loc\big([0,\infty);Y^\star\big).$ Finally, by the dense embedding $\Dr(\Omega)\inj Y,$ we deduce that $L^1_\loc\big([0,\infty);Y^\star\big)\inj\Dr^\p\big((0,\infty)\times\Omega\big)$ and the proposition is proved.
\medskip
\end{vproof}

\begin{vproof}{of Proposition~\ref{propdep}.}
As we shall see, the proof can be easily from the one given in~Bégout and D\'iaz~\cite[Lemma~6.5]{MR4053613}. The embedding in \eqref{propdephyp} comes from \eqref{Xcon}. We make the difference between the two equations satisfied by $u$ and $\wt u.$ If follows from Lemmas~\ref{lemVL} and \ref{lemcong} that $u-\wt u$ satisfies the equation obtained in $L^1_\loc\big((0,\infty);X^\star\big).$ We take the $X^\star-X$ duality product with $\vi(u-\wt u).$ By Corollaries~\ref{corlemineg}, \ref{corA0mon}, (A.3) of Lemma~A.5 in~Bégout and D\'iaz~\cite{MR4053613} and Cauchy-Schwarz's inequality, we then arrive at,
\begin{gather*}
\frac12\frac\d{\d t}\|u-\wt u\|_{L^2(\Omega)}^2\le\|f-\wt f\|_{L^2(\Omega)}\|u-\wt u\|_{L^2(\Omega)},
\end{gather*}
almost everywhere on $(0,\infty).$ Integrating over $(s,t),$ we obtain \eqref{estthmweak}.
\medskip
\end{vproof}

\begin{vproof}{of Theorem~\ref{thmstrongH2}.}
Let the assumptions of the theorem be fulfilled. By Corollary~\ref{corAmaxmon} and Barbu~\cite[Theorem~2.2, p.131]{MR0390843} (see also Vrabie~\cite[Theorem~1.7.1, p.23]{MR1375237}), there exists a unique $u\in W^{1,\infty}_\loc\big([0,\infty);L^2(\Omega)\big)$ satisfying $u(t)\in D(A_0^m)$ and \eqref{nls} in $L^2(\Omega),$ for almost every $t>0,$ $u(0)=u_0$ and \eqref{strongH23}. This last estimate yields \eqref{strongH21}. Since $u\in W^{1,\infty}_\loc\big([0,\infty);L^2(\Omega)\big),$ it follows from Lemma~A.5 in~Bégout and D\'iaz~\cite{MR4053613} that the map $M:t\longmapsto\frac12\|u(t)\|_{L^2(\Omega)}^2$ belongs to $W^{1,\infty}_\loc\big([0,\infty);\R\big)$ and $M^\p(t)=\big(u(t),u_t(t)\big)_{L^2(\Omega)},$ for almost every $t>0.$ Taking the $L^2$-scalar product of \eqref{nls} with $\vi u,$ we obtain \eqref{L2}, for almost every $t>0.$ By \eqref{L2} and Cauchy-Schwarz's inequality, we get
\begin{gather}
\label{demH21}
\forall t\ge0, \; \|u(t)\|_{L^2(\Omega)}\le\|u_0\|_{L^2(\Omega)}+\int_0^t\|f(s)\|_{L^2(\Omega)}\d s,
\end{gather}
We multiply \eqref{L2} by $C_0=\frac{|\Re(a)|+1}{\Im(a)}.$ Then, we take again the $L^2$-scalar product of \eqref{nls} with $u.$ Summing the result with $C_0\times\eqref{L2},$ we infer
\begin{gather*}
\|\nabla u\|_{L^2(\Omega)}^2+\|u\|_{L^{m+1}(\Omega)}^{m+1}
\le C\left(\|u_t\|_{L^2(\Omega)}+\|Vu\|_{L^2(\Omega)}+\|f\|_{L^2(\Omega)}\right)\|u\|_{L^2(\Omega)},
\end{gather*}
almost everywhere on $(0,\infty).$ It follows from \eqref{lemV1}--\eqref{lemV2} that for some $C=C(N,C_0)$ $(C=C(\beta,C_0),$ if $N=2),$
\begin{gather}
\label{demH22}
\begin{split}
	&	\; \|\nabla u\|_{L^2(\Omega)}^2+\|u\|_{L^{m+1}(\Omega)}^{m+1}			\\
  \le	&	\; C\left(\|u_t\|_{L^2(\Omega)}+\left(\|V_1\|_{L^\infty(\Omega)}+\|V_2\|_{L^{p_V}(\Omega)}^{2-\gamma}\right)\|u\|_{L^2(\Omega)}
			+\|f\|_{L^2(\Omega)}\right)\|u\|_{L^2(\Omega)},
\end{split}
\end{gather}
almost everywhere on $(0,\infty).$ By \eqref{demH21}--\eqref{demH22}, $u\in L^\infty_\loc\big([0,\infty);H^1_0(\Omega)\cap L^{m+1}(\Omega)\big)$ and $u$ is an $H^2$-solution. Using \eqref{lemVL2}, we get $Vu\in L^\infty_\loc\big([0,\infty);L^2(\Omega)\big).$ By \eqref{nls}, if $m=\{1,0\}$ then $\Delta u\in L^\infty_\loc\big([0,\infty);L^2(\Omega)\big).$ Now, asssume that $0<m<1.$ Since $a\in C_{\mathrm{int}}(m),$ there exists $b\in\C$ such that $|b|=1,$ $\Im(b)<0$ and $ab\in C_{\mathrm{int}}(m)$ (Bégout~\cite[Lemma~4.2]{MR4098330}). We take the $L^2$-scalar product of \eqref{nls} with $\vi abg_0^m(u).$ We get,
\begin{equation*}
\begin{aligned}
	&	\; \Re\left(\vi ab\vint_{\Omega}\ovl{\vi u_t}g_0^m(u)\d x\right)+\Re\left(\vi ab\vint_{\Omega}g_0^m(u)\ovl{\Delta u}\d x\right)
			+\Re(\vi ab)\vint_\Omega V\ovl ug_0^m(u)\d x		\\
   +	&	\; |a|^2|\Im(b)|\|g_0^m(u)\|_{L^2(\Omega)}^2=\Re\left(\vi ab\vint_{\Omega}\ovl fg_0^m(u)\d x\right).
\end{aligned}
\end{equation*}
By Lemma~\ref{lemenereg}, the second term in the left hand side of the above is nonnegative which becomes,
\begin{gather*}
|a||\Im(b)|\,\|u\|_{L^{2m}(\Omega)}^{2m}\le\vint_{\Omega}|\vi u_t+Vu-f|\,|g(u)|\d x.
\end{gather*}
By Cauchy-Schwarz's and Young's inequalities, we get
\begin{gather}
\label{demH23}
\|u\|_{L^{2m}(\Omega)}^{2m}\le\frac1{(|a||\Im(b)|)^2}\|\vi u_t+Vu-f\|_{L^2(\Omega)}^2,
\end{gather}
almost everywhere on $(0,\infty).$ We deduce from \eqref{nls}and \eqref{demH23} that,
\begin{gather}
\label{demH24}
\|\Delta u\|_{L^2(\Omega)}+|a|\|u\|_{L^{2m}(\Omega)}^m\le\frac2{|\Im(b)|}\|\vi u_t+Vu-f\|_{L^2(\Omega)},
\end{gather}
almost everywhere on $(0,\infty).$ It follows that $u\in L^\infty_\loc\big([0,\infty);L^{2m}(\Omega)\big)$ and $\Delta u\in L^\infty_\loc\big([0,\infty);L^2(\Omega)\big).$ Now we go back to the general case $0\le m\le1.$ Then \eqref{strongH22} follows from \eqref{strongH21} and the estimate,
\begin{gather}
\label{nablau}
\|\nabla u\|_{L^2(\Omega)}^2\le\|u\|_{L^2(\Omega)}\|\Delta u\|_{L^2(\Omega)},
\end{gather}
which holds for any $u\in H^1_0(\Omega)$ such that $\Delta u\in L^2(\Omega).$ Finally, the rest of the properties is clear by \eqref{nls}, \eqref{strongH23}, \eqref{lemVL2},  \eqref{demH21}, \eqref{demH22}, \eqref{demH24} and Remark~\ref{rmkf0}. This ends the proof of the theorem.
\medskip
\end{vproof}

\begin{vproof}{of Theorem~\ref{thmweak}.}
Existence, estimate~\eqref{estthmweak} and uniqueness come from density of $\Dr(\Omega)\times W^{1,1}_\loc([0,\infty);L^2(\Omega))$ in $L^2(\Omega)\times L^1_\loc([0,\infty);L^2(\Omega)),$ Theorem~\ref{thmstrongH2}, Proposition~\ref{propdep} and completeness of $C\big([0,T];L^2(\Omega)\big),$ for any $T>0.$ Let $u$ be the unique weak solution. Then $u$ is a limit of $H^2$-strong solutions $(u_n)_{n\in\N}$ in $C([0,T];L^2(\Omega)),$ for any $T>0.$ By \eqref{L2}, each $u_n$ satisfies \eqref{L2+} with equality. If $|\Omega|<\infty$ or if $m=1$ then we can pass to the limit to obtain \eqref{Lm}--\eqref{L2+}, still with equality. Otherwise, we work with $\|u_n(\sigma)\|_{L^{m+1}(\Omega\cap B(0,R))}^{m+1}$ in place of $\|u_n(\sigma)\|_{L^{m+1}(\Omega)}^{m+1}$ in \eqref{L2+}, pass to the limit in $n$ and then in $R.$ For more details, see the proof of Bégout and D\'{\i}az~\cite[Proposition~2.3]{MR4053613}.
\medskip
\end{vproof}

\begin{vproof}{of Theorems~\ref{thmstrongaH1} and \ref{thmstrongH1}.}
Let $u_0\in H^1_0(\Omega)$ and $f\in L^1_\loc\big([0,\infty);H^1_0(\Omega)\big).$ Let $(\vphi_n)_{n\in\N}\subset\Dr(\Omega)$ and $(f_n)_{n\in\N}\subset\Dr\big([0,\infty);H^1_0(\Omega)\big)$ be such that $\vphi_n\xrightarrow[n\to\infty]{H^1_0(\Omega)}u_0$ and $f_n\xrightarrow[n\to\infty]{L^1((0,T);H^1_0)}f,$ for any $T>0.$ For each $n\in\N,$ let $u_n$ be the unique $H^2$-solution to \eqref{nls} such that $u_n(0)=\vphi_n,$ given by Theorem~\ref{thmstrongH2}. By Proposition~\ref{propdep}, $(u_n)_{n\in\N}$ is a Cauchy sequence in $C\big([0,T];L^2(\Omega)\big),$ for any $T>0.$ As a consequence, there exists $u\in C\big([0,\infty);L^2(\Omega)\big)$ such that for any $T>0,$
\begin{gather}
\label{demthmstrongH12}
u_n\xrightarrow[n\to\infty]{C([0,T];L^2(\Omega))}u.
\end{gather}
By definition, $u$ is a weak solution and satisfies \eqref{nls} in $\Dr\big((0,\infty)\times\Omega\big)$ (Proposition~\ref{propsolL2}). In particular, $u$ fulfills \eqref{Lm}. Taking the $L^2$-scalar product of \eqref{nls} with $-\vi\Delta u_n,$ it follows from~\cite[Lemma~A.5]{MR4053613} and Lemma~\ref{lemenereg} that for any $n\in\N$ and almost every $s>0,$
\begin{gather*}
\frac12\frac{\d}{\d t}\|\nabla u_n(s)\|_{L^2(\Omega)}^2\le\big(\nabla f_n(s)-u_n(s)\nabla V,\vi\nabla u_n(s)\big)_{L^2(\Omega)}.
\end{gather*}
We apply Cauchy-Schwarz's inequality, \eqref{lemVL2} and \eqref{L2}. We get for any $n\in\N$ and almost every $s>0,$
\begin{gather*}
\frac12\frac{\d}{\d t}\|u_n(t)\|_{H^1_0}^2\le\|f_n(t)\|_{H^1_0}\|u_n(t)\|_{H^1_0}+C\|\nabla V\|_{L^\infty+L^{p_V}}\|u_n(t)\|_{H^1_0}^2,
\end{gather*}
where $C$ is given by \eqref{lemVL2}. After integration, we obtain
\begin{gather*}
\|u_n(t)\|_{H^1_0}\le\|\vphi_n\|_{H^1_0}+\vint_0^t\|f_n(s)\|_{H^1_0}\d s+\vint_0^tC\|\nabla V\|_{L^\infty+L^{p_V}}\|u_n(s)\|_{H^1_0}\d s,
\end{gather*}
and by Gronwall's Lemma,
\begin{gather}
\label{demthmstrongH13}
\|u_n(t)\|_{H^1_0}\le\left(\|\vphi_n\|_{H^1_0}+\vint_0^t\|f_n(s)\|_{H^1_0}\d s\right)e^{C\|\nabla V\|_{L^\infty+L^{p_V}}t},
\end{gather}
for almost every $t>0$ and any $n\in\N.$ It follows that,
\begin{gather}
\label{demthmstrongH14}
(u_n)_{n\in\N} \text{ is bounded in } C\big([0,T];H^1_0(\Omega)\big),
\end{gather}
for any $T>0.$ By \eqref{demthmstrongH14}, \eqref{lemVL2} and \eqref{L2+},
\begin{gather}
\label{demthmstrongH15}
\big(\Delta u_n+Vu_n+ag(u_n)\big)_{n\in\N}  \text{ is bounded in } L^\frac{m+1}m_\loc\big([0,\infty);X^\star\big),
\end{gather}
for any $T>0,$ where $X^\star=H^{-1}(\Omega)+L^\frac{m+1}m(\Omega).$ We have $L^2\big((0,T);L^2(\Omega)\big)\inj L^1\big((0,T);H^{-1}(\Omega)\big)$ with dense embedding and $L^2\big((0,T);L^2(\Omega)\big)\cong L^2((0,T)\times\Omega),$ which is separable. It follows that $L^1\big((0,T);H^{-1}(\Omega)\big)$ is separable, for any $T>0.$ In addition, $H^{-1}(\Omega)$ is a reflexive Banach space, and so $L^1\big((0,T);H^{-1}(\Omega)\big)^\star\cong L^\infty\big((0,T);H^1_0(\Omega)\big),$ for any $T>0$ (Edwards~\cite[Theorem~8.18.3, p.590]{MR0221256}). With help of \eqref{demthmstrongH12} and \eqref{demthmstrongH14}, it follows that $u\in L^\infty_\loc\big([0,\infty);H^1_0(\Omega)\big)$ and for any $T>0,$
\begin{gather}
\label{demthmstrongH16}
u_n\underset{n\to\infty}{-\!\!\!-\!\!\!-\!\!\!-\!\!\!\weak}u, \text{ in } L^\infty_{\w\star}\big((0,T);H^1_0(\Omega)\big).
\end{gather}
We deduce from \eqref{demthmstrongH12}, \eqref{demthmstrongH16}, \eqref{lemVL2}, \eqref{Lm} and Lemma~\ref{lemcong} that $u$ satisfies the first line of \eqref{thmstrongaH11} and
\begin{gather}
\label{demthmstrongH17}
\Delta u+Vu+ag(u)\in L^\frac{m+1}m_\loc\big([0,\infty);X^\star\big).
\end{gather}
By \eqref{demthmstrongH17} and \eqref{nls}, $u$ satisfies \eqref{thmstrongaH11}, and by \eqref{demthmstrongH16}, \eqref{demthmstrongH13} and the lower semicontinuity of the norm, $u$ satisfies \eqref{thmstrongaH12} with $s=0.$ Now, we fix $s>0.$ Let $v$ be the unique weak solution to \eqref{nls} given by this proof, where $v(0)=u(s)$ and $t\longmapsto f(t)$ is replaced with $t\longmapsto f(t+s).$ By uniqueness, $v(t)=u(t+s),$ for any $t\ge0.$ We then obtain the general case \eqref{thmstrongaH12}. Finally, the continuous embedding in $\Dr^\p\big((0,\infty)\times\Omega\big)$ stated in Theorem~\ref{thmstrongaH1} comes from the dense embedding $\Dr(\Omega)\inj H^1_0(\Omega)\cap L^{m+1}(\Omega).$ This ends the proof of Theorem~\ref{thmstrongaH1}. \\
Now assume further that $f\in L^\frac{m+1}m_\loc\big([0,\infty);X^\star\big).$ It follows from \eqref{nls}, \eqref{demthmstrongH15}  and \eqref{demthmstrongH17} that
\begin{gather}
\nonumber
u\in W^{1,\frac{m+1}m}_\loc\big([0,\infty);X^\star\big),					\\
\label{demthmstrongH18}
(u_n)_{n\in\N} \text{ is bounded in } W^{1,\frac{m+1}m}\big((0,T);X^\star\big),
\end{gather}
for any $T>0.$ Hence $u$ is an $H^1_0$-solution. Using the embedding,
\begin{gather*}
W^{1,\frac{m+1}m}\big((0,T);X^\star\big)\inj C^{0,\frac1{m+1}}\big([0,T];X^\star\big),
\end{gather*}
it follows from \eqref{demthmstrongH12}, \eqref{demthmstrongH14}, \eqref{demthmstrongH18} and Cazenave~\cite[Proposition~1.1.2, p.2]{MR2002047} that,
\begin{gather*}
u\in C_\w\left([0,\infty);H^1_0(\Omega)\right).
\end{gather*}
By \ref{rmkdefsol2} of Remark \ref{rmkdefsol}, we can take the $X-X^\star$ duality product of \eqref{nls} with $\vi u.$ Applying~\cite[Lemma~A.5]{MR4053613}, Property~\ref{thmstrongH12} follows. Finally, $u$ is the unique $H^1_0$-solution by Proposition~\ref{propdep} (and also by \eqref{rmkdefsol11} if $m=0).$ This ends the proof of Theorem~\ref{thmstrongH1}.
\medskip
\end{vproof}

\section{Proofs of the finite time extinction and asymptotic behavior theorems}
\label{proofext}

In this section, we shall prove the results of Section~\ref{finite}. Theorems~\ref{thmextN1}, \ref{thmrtdH1}, \ref{thmextH2} and \ref{thmrtdH2} may be obtained with the same method, while Theorems~\ref{thmextN1e} and \ref{thmextH2e} require an adaptation. As far as we know, the pionering result to obtain finite time extinction for solutions of some damped nonlinear Schrödinger equation is due to Carles and Gallo~\cite{MR2765425}. As said in the Introduction, the present extension is possible thanks to a sharper study of the regularity and existence frameworks. In addition, synchronized finite extinction time and the results for $m=0$ and $f(t)$ non zero are completely new.
\medskip

\begin{vproof}{of Theorems~\ref{thmextN1}--\ref{thmrtdH1}, \ref{thmextH2} and \ref{thmextH2e}--\ref{thmrtdH2}.}
The proof of these Theorems relies on the following Gagliardo-Nirenberg inequality which asserts that there exists $C_\GN=C_\GN(m,N)$ such that for any $u\in H^1_0(\Omega)\cap L^{m+1}(\Omega),$
\begin{gather}
\label{GN1}
\|u\|_{L^2(\Omega)}^\frac{(N+2)-m(N-2)}2\le C_\GN\|u\|_{L^{m+1}(\Omega)}^{m+1}\|\nabla u\|_{L^2(\Omega)}^\frac{N(1-m)}2.
\end{gather}
If, in addition, $\Delta u\in L^2(\Omega),$ then it follows from \eqref{nablau} that
\begin{gather}
\label{GN2}
\|u\|_{L^2(\Omega)}^\frac{(N+4)-m(N-4)}4\le C_\GN\|u\|_{L^{m+1}(\Omega)}^{m+1}\|\Delta u\|_{L^2(\Omega)}^\frac{N(1-m)}4.
\end{gather}
Now, suppose Assumptions~\ref{assN1} or \ref{assH2} are fulfilled. In Theorems~\ref{thmextN1} and \ref{thmrtdH1}, by \eqref{f}, $u$ becomes a strong solution (except for $m=0).$ Therefore, \eqref{L2} is satisfied on $(T_0,\infty)$ (which comes from the equality \eqref{L2+}, if $m=0).$ In Theorems~\ref{thmextN1e}, \ref{thmextH2}, \ref{thmextH2e} and \ref{thmrtdH2}, $u$ is always a strong solution and \eqref{L2} is verified almost everywhere on $(0,\infty)$ Now, we let $\ell=1$ for the proof of Theorems~\ref{thmextN1}--\ref{thmrtdH1}, and $\ell=2$ for the proof of Theorems~\ref{thmextH2} and \ref{thmextH2e}--\ref{thmrtdH2}. By \eqref{thmstrongaH12} and Theorem~\ref{thmstrongH2}, it follows that $u\in L^\infty\big((0,\infty);H^1_0(\Omega)\big),$ if $\ell=1,$ with additionally $\Delta u\in L^\infty\big((0,\infty);L^2(\Omega)\big),$ if $\ell=2.$ Setting for any $t\ge0,$ $y(t)=\|u(t)\|_{L^2(\Omega)}^2,$ $\alpha_1=\Im(a)C_\GN^{-1}\|\nabla u\|_{L^\infty((0,\infty);L^2(\Omega))}^{-\frac{N(1-m)}2},$ and $\alpha_2=\Im(a)C_\GN^{-1}\|\Delta u\|_{L^\infty((0,\infty);L^2(\Omega))}^{-\frac{N(1-m)}4},$ it follows from \eqref{L2}, \eqref{GN1}--\eqref{GN2} and Cauchy-Schwarz's inequality that,
\begin{gather}
\label{edo}
y^\p(t)+2\alpha_\ell y(t)^{\delta_\ell}\le2\|f(t)\|_{L^2(\Omega)}y(t)^\frac12,
\end{gather}
where $\delta_\ell$ is defined by \eqref{delta}. We proceed with the proof in four steps.
\\
\textbf{Step~1:} Proof of Theorem~\ref{thmextN1}, \ref{thmrtdH1}, \ref{thmextH2} and \ref{thmrtdH2}. \\
By Assumptions~\ref{assN1} and \ref{assH2}, if $m\neq0$ then the right hand side member of \eqref{edo} vanishies on $(T_0,\infty)$ and after integration, we obtain the results of these theorems. If $m=0,$ we have by \eqref{L2} and Hölder's inequality,
\begin{gather*}
\frac12\frac\d{\d t}\|u(t)\|_{L^2(\Omega)}^2+\omega_f\|u(t)\|_{L^1(\Omega)}\le0,
\end{gather*}
for almost every $t>T_0,$ with $\omega_f=\Im(a)-\|f\|_{L^\infty((T_0,\infty)\times\Omega)}.$ From assumption \eqref{f} we know that $\omega_f>0$. Then, by
the Gagliardo-Nirenberg interpolation inequalities \eqref{GN1}--\eqref{GN2} we get that for almost every $t>T_0,$
\begin{gather*}
y^\p(t)+\beta_\ell y(t)^{\delta_\ell}\le0,
\end{gather*}
where $\beta_1=2\omega_fC_\GN^{-1}\|\nabla u\|_{L^\infty((0,\infty);L^2(\Omega))}^{-\frac{N}2},$ and $\beta_2=2\omega_fC_\GN^{-1}\|\Delta u\|_{L^\infty((0,\infty);L^2(\Omega))}^{-\frac{N}4}.$ And again the conclusion follows by integration.
\medskip \\
We turn out to the proof of Theorems~\ref{thmextN1e} and \ref{thmextH2e}. Let $\alpha=\Im(a)C_\GN^{-1}.$
\\
\textbf{Step~2:} In Theorem~\ref{thmextN1e}, there exists $\eps_\star=\eps_\star(|a|,m)$ with,
\begin{gather}
\label{eps*}
\eps_\star\le\min\left\{(2\delta_\ell-1)^{-\frac{2\delta_\ell-1}{\delta_\ell}}(\alpha\delta_\ell)^\frac1{1-\delta_\ell}(1-\delta_\ell)^\frac{2\delta_\ell-1}{\delta_\ell(1-\delta_\ell)},\alpha\,\delta_\ell\,(1-\delta_\ell)\right\},
\end{gather}
such that if \eqref{thmextN1e1} holds true then $\|\nabla u\|_{L^\infty((0,\infty);L^2(\Omega))}\le1.$
\\
This comes from \eqref{thmstrongaH12V0}.
\\
\textbf{Step~3:} In Theorem~\ref{thmextH2e}, there exists $\eps_\star=\eps_\star(|a|,N,m)$ satisfying \eqref{eps*} such that under assumption \eqref{thmextH2e1}, we have $\|\Delta u\|_{L^\infty((0,\infty);L^2(\Omega))}\le1.$ \\
This comes from \eqref{estthmweak}, \eqref{strongH23}, \eqref{lemVL2}, \eqref{demH22} and \eqref{demH24}.
\\
\textbf{Step~4:} Proof of Theorems~\ref{thmextN1e} and \ref{thmextH2e}.
\\
Let $x_\star=(\alpha\delta_\ell(1-\delta_\ell)T_0)^\frac1{1-\delta_\ell}$ and $y_\star=\big(\alpha\delta_\ell^{\delta_\ell}(1-\delta_\ell)\big)^\frac1{1-\delta_\ell}.$ By \eqref{thmextN1e1}, \eqref{thmextH2e1} and \eqref{eps*},
\begin{gather}
\label{demthmextH2e1}
y(0)\le x_\star.
\end{gather}
By Steps~2 and 3, $\alpha\le\min\big\{\alpha_1,\alpha_2\big\}.$ Applying Young's inequality to \eqref{edo}, we arrive at,
\begin{gather*}
\begin{split}
	&	\; y^\p(t)+2\alpha y(t)^{\delta_\ell}											\medskip \\
  \le	&	\; \frac{2\delta_\ell-1}{\delta_\ell}(\alpha\delta_\ell)^{-\frac1{2\delta_\ell-1}}
			\|f(t)\|_{L^2(\Omega)}^\frac{2\delta_\ell}{2\delta_\ell-1}+\alpha y(t)^{\delta_\ell},
\end{split}
\end{gather*}
for almost every $t\ge0.$ Replacing \eqref{thmextN1e1} and \eqref{thmextH2e1} in the above and using \eqref{eps*}, we obtain
\begin{gather}
\label{demthmextH2e2}
y^\p(t)+\alpha y(t)^{\delta_\ell}\le y_\star\big(T_0-t\big)_+^\frac{\delta_\ell}{1-\delta_\ell},
\end{gather}
for almost every $t>0.$ By \eqref{demthmextH2e1}, \eqref{demthmextH2e2} and \cite[Lemma~5.2]{MR4053613}, $y(t)=0,$ for any $t\ge T_0.$
\medskip
\end{vproof}

\begin{vproof}{of Theorem~\ref{thm0w}.}
By \eqref{estthmweak}, density and Remark~\ref{rmkthm0w}, we may assume that $f\in\Dr\big([0,\infty);L^2(\Omega)\big),$ $u_0\in\Dr(\Omega)$ and $m<1.$ The result then comes from Theorems~\ref{thmextH2} and \ref{thmrtdH2}.
\medskip
\end{vproof}

\begin{vproof}{of Theorem~\ref{thmtdH1}.}
By \eqref{thmstrongaH12}, $u\in L^\infty\big((0,\infty);H^1_0(\Omega)\big).$ The result then comes from Theorem~\ref{thm0w} and Gagliardo-Nirenberg's inequality.
\medskip
\end{vproof}

\begin{vproof}{of Theorem~\ref{thmtdH2}.}
By Property~\ref{thmstrongH24} of Theorem~\ref{thmstrongH2}, Theorem~\ref{thm0w} and \eqref{nablau}, $\vlim_{t\nearrow\infty}\|u(t)\|_{H^1_0(\Omega)}=0.$ The second limit is due to the first one, Hölder's inequality and the Sobolev embeddings. The last limit comes from the two first and \eqref{L2}.
\medskip
\end{vproof}

\section*{Acknowledgements}
\baselineskip .5cm
P.~Bégout acknowledges funding from ANR under grant ANR-17-EURE-0010 (Investissements d’Avenir program). The research of J.~I.~D\'{\i}az was partially supported by the projects ref. MTM2017-85449-P and PID2020-112517GB-I00 of the DGISPI (Spain) and the Research Group MOMAT (Ref. 910480) of the UCM.

\baselineskip .4cm

\def\cprime{$^\prime$}

\end{document}